\documentclass[11pt]{article}
\usepackage{amssymb,latexsym,amsfonts,verbatim,amscd}
\usepackage[all]{xy}

\newtheorem{Def}{Definition}[section]
\newtheorem{Thm}[Def]{Theorem}
\newtheorem{Lem}[Def]{Lemma}
\newtheorem{Prop}[Def]{Proposition}
\newtheorem{Cor}[Def]{Corollary}

\newtheorem{Fac}[Def]{Fact}

\def\telos{\hfill$\dashv$}

\begin{document}

\title{Semantics for first-order superposition logic}
\author{Athanassios Tzouvaras}

\date{}
\maketitle

\begin{center}
 Department  of Mathematics\\  Aristotle University of Thessaloniki \\
 541 24 Thessaloniki, Greece \\
  e-mail: \verb"tzouvara@math.auth.gr"
\end{center}

\begin{abstract}
We investigate how the  sentence choice semantics (SCS) for propositional superposition logic (PLS) developed in \cite{Tz17} could be extended so as to successfully apply to first-order superposition logic(FOLS). There are two options for  such an  extension. The apparently more natural one is  the formula choice  semantics (FCS) based on choice functions for pairs of arbitrary formulas of the basis language. It is proved however that the universal instantiation scheme of FOL, $(\forall v)\varphi(v)\rightarrow\varphi(t)$, is false, as a scheme of tautologies, with respect to FCS. This causes  the total failure of FCS as a candidate semantics.  Then we turn to the other option  which is a variant of SCS, since it uses again choice functions for pairs of sentences only. This semantics however presupposes  that the  applicability of the connective $|$ is restricted to quantifier-free sentences, and thus the class of well-formed formulas and sentences of the language is restricted too. Granted  these syntactic restrictions, the usual axiomatizations of FOLS turn out to be  sound and conditionally complete with respect to this second semantics, just like the corresponding systems of PLS.
\end{abstract}

\vskip 0.2in

{\em Mathematics Subject Classification (2010)}: 03B60, 03G12

\vskip 0.2in

{\em Keywords:}  Logic of superposition.  Choice function for pairs of sentences/formulas. Sentence choice semantics. Formula choice semantics.

\section{Introduction}
In \cite{Tz17} we introduced and investigated various systems of propositional superposition logic (PLS). The systems of PLS extend classical propositional logic (PL). Their language  is that of PL augmented with a new binary connective $|$ for the ``superposition operation'', while their axioms are those of PL together with a few axioms about $|$.  The motivating idea was roughly this: if $\varphi|\psi$ denotes the ``superposition of two states'' (or, more precisely, the propositions expressing  these states), as the latter is currently understood  in quantum mechanics (QM), what is the  {\em purely logical content} of the operation, that is, what can we say about the {\em truth} of $\varphi|\psi$ without leaving  the ground of classical logic? The basic intuition  is that $\varphi|\psi$ strangely expresses  {\em both} some kind of   conjunction of the  properties $\varphi$ and $\psi$ (before the measurement), and simultaneously  some kind of   disjunction of the same properties (after the measurement, i.e., after the ``collapsing'' of the superposed states). This collapsing can be formalized by the help of a choice function that acts on pairs of sentences $\{\varphi,\psi\}$, turning each formula $\varphi|\psi$ into a classical one. Such functions formed the basis of a semantics  for the new logic, called {\em sentence choice semantics} (or SCS for short), that allows  $\varphi|\psi$ to present simultaneously conjunctive and disjunctive characteristics, which are manifested in the ``interpolation property'', i.e., the property of  $\varphi|\psi$ to be strictly logically interpolated between $\varphi\wedge\psi$ and $\varphi\vee\psi$.

Although QM has been the source of motivation for introducing the  logical connective of superposition, PLS is {\em not} a quantum logic based on orthomodular lattices (see \cite{Ka83} for complete information about such  structures), as these logics are  discussed e.g. in  \cite{PM09}. Nor is it in a similar  vein with the content, e.g., of  \cite{BS04}, \cite{Ba10}, and other papers cited and discussed in \cite{Tz18}, that belong to what can be called standard approach to quantum phenomena.  As  said above, the logic(s) PLS  intend only  to explore the logical content of the phenomenon of superposition  alone. To quote from \cite[p. 151]{Tz17}:
\begin{quote}
``So  the logic presented here is hardly  the logic of superposition as this concept is currently used and understood in  physics today. It is rather the logic of superposition, when the latter is  understood as the `logical extract' of the corresponding physics concept. Whether it could eventually have applications to the field of QM we don't know.''
\end{quote}
In response to questions asked by one of the referees let me add some further comments. That PLS (and its  first-order extension FOLS considered in this paper) has almost no points in common with standard treatments of QM can be  simply inferred from the fact that  it is {\em not} a probabilistic theory. Probabilities have no place in this logical system  (as it stands) and I cannot see how it  could be revised in order to be compatible with their use. This is why the ``collapsing'' of the superposed sentence $\varphi|\psi$ is accomplished by means of a {\em choice} between $\varphi$ and $\psi$. And up to my knowledge there is no genuine theory that relates fruitfully choice functions with probabilities. As we put it in \cite[p. 151]{Tz17}:
\begin{quote}
``The basic idea is that the collapse  of the composite state  $c_0\vec{u}_0+c_1\vec{u}_1$ to one of the states $\vec{u}_0$, $\vec{u}_1$  can be seen, {\em from the point of view of pure logic,} just as a (more or less random)  choice from  the set of possible outcomes $\{\vec{u}_0,\vec{u}_1\}$. This is because from the point of view of pure logic   probabilities are  irrelevant or, which amounts to the same thing, the states $\vec{u}_0$ and $\vec{u}_1$ are considered  equiprobable.   In such a case  the  superposition of $\vec{u}_0$ and $\vec{u}_1$ is unique and the outcome of the collapse can be decided  by a coin tossing or, more strictly,  by  a  {\em choice function} acting on  pairs of observable states, which in our case coincide with pairs of sentences of  $L$. This  of course constitutes  a major deviation from the standard treatment of superposition, according to which there is not just one superposition of  $\vec{u}_0$ and  $\vec{u}_1$ but infinitely many, actually as many  as the number of  linear combinations $c_0\vec{u}_0+c_1\vec{u}_1$, for $|c_0|^2+|c_1|^2=1$.''
\end{quote}
Of course, theoretically,  we could switch  from $\{0,1\}$-valuations of classical logic to $[0,1]$-valuations of a non-classical logic. But then the interpretation of superposition would not be ``within classical reasoning and commonsense'', as was the aim of the original attempt. Perhaps in the future we shall attempt some non-classical interpretation through a continuous-valued logic.

In view of the above, fundamental concepts pertaining to the probabilistic character of the standard treatment of QM, such as global vs local phases, contextuality, the Kochen-Specker theorem, entanglement, etc., simply do not  make any sense for our logic. However, despite of this,  the systems PLS still seem to have merits. As a second  referee wrote: ``Even though the interpretation of superposition logics in terms of the original quantum-mechanical motivation is probably dubious, still the connective with a choice semantics is sufficiently interesting per se to justify the investigation; conceivably, the logic can find other interpretations, perhaps of some epistemic or possibilistic kind.''

Let us come now to the content of the present work. A natural question, already asked in the last section of \cite{Tz17}, is whether  the logic of superposition can be extended to a quantified version, i.e., whether the systems of PLS can be extended to corresponding systems of first-order superposition logic (or FOLS for short). At a syntactic level, systems of FOLS extending corresponding systems of PLS are very easily defined. They are just extensions of classical first-order logic (FOL) by the help of the same axioms for $|$ that were used in PLS. This is because there are no new axioms for $|$ involving $\forall$ or $\exists$, as  there are no plausible correlations between $|$ and the quantifiers.
But at semantic level things are much more complex. First, we made sure that systems of FOLS do have semantics having characteristics quite analogous to that of SCS for PLS. Actually, an alternative semantics that we meanwhile developed  for PLS in \cite{Tz18}, the Boolean-value choice semantics (or BCS for short), turned out to be suitable also for FOLS.

However the question whether a semantics for FOLS generalizing  SCS of PLS is possible, remained. One of the goals  of the present paper is to show that the straightforward generalization of the semantics SCS of PLS, namely the semantics based on choice functions for all pairs of {\em formulas} of a first-order language $L$, called {\em formula choice semantics} (or FCS for short), does not work. Specifically, we show that the systems of FOLS fail to be true with respect to FCS (i.e., soundness fails) in an unexpected way: It is not the axioms for $|$ that fail to be tautologies of FCS but one of the fundamental axiom of FOL, namely the Universal Instantiation ($UI$) scheme, $\forall v\varphi(v)\rightarrow \varphi(t)$. This of course leads to the break down of FCS itself, since it cannot accommodate  the most fundamental logical constant of quantified logic, the universal quantifier. This result is shown in section 3. In section 4 we show a related fact concerning non-existence of uniform choice functions.

The second major result shown in this paper  is that the semantics SCS (using functions on pairs of {\em sentences} only rather than arbitrary formulas) {\em can} be applied also to FOLS, provided we shall restrict the applicability of the connective $|$ to formulas without quantifiers (unless they are classical ones), and thus restrict the class of well-formed formulas of the language $L_s=L\cup\{|\}$ of the logic of superposition. Under this restriction, FOLS is proved sound and conditionally complete with respect to SCS. This result is described in section 5.

Since the content of the present paper relies heavily on the material contained in \cite{Tz17}, we need  first to  recall briefly  the main notions and facts established there. This is done in the next subsection.

\subsection{Overview of PLS with sentence choice semantics}
This  subsection overviews the main notions and facts contained  in \cite{Tz17}. It is identical to the corresponding introductory subsection 1.1 of \cite{Tz18}. In general a Propositional Superposition Logic (PLS) consists, roughly, of a pair $(X,K)$, where  $X$ is the  semantic and $K$ is the syntactic  part of the logic.  Actually $K$ is a {\em formal system} in the usual sense of the word, and $X$ is a set of functions that provides meaning to sentences in a way described below.   ${\rm PLS}(X,K)$ will denote the propositional superposition logic with semantic part $X$ and syntactic part $K$. The precise  definition of ${\rm PLS}(X,K)$ will be given  below.

Although the semantic part is the most intuitively appealing we start with  the description of the syntactic part $K$. The  language of all formal systems  $K$ below (or the language of PLS), $L_s$, is that  of standard Propositional Logic (PL) $L=\{p_0,p_1,\ldots\}\cup\{\wedge,\vee,\rightarrow,\leftrightarrow,\neg\}$ augmented with the new binary connective ``$|$''. That is,  $L_s=L\cup\{|\}$.
The set of sentences of $L_s$, $Sen(L_s)$, is defined by induction as usual, with the additional inductive step that $\varphi|\psi$ is a sentence whenever $\varphi$ and $\psi$ are so.

Throughout the letters   $\alpha$, $\beta$, $\gamma$ range exclusively over the set of sentences of $L$, $Sen(L)$,  while $\varphi$, $\psi$, $\sigma$  range over elements of $Sen(L_s)$ in general.

A formal system  $K$ consists of a set of axioms   $\textsf{Ax}(K)$ and a set of inference rules $\textsf{IR}(K)$.  The axioms of $K$ always include the axioms of PL, while $\textsf{IR}(K)$ includes the inference rule of PL. So let us first fix the axiomatization for PL consisting of the following axiom schemes (for the language $L_s$).

\vskip 0.1in

(P1) $\varphi\rightarrow (\psi\rightarrow\varphi)$

(P2) $(\varphi\rightarrow(\psi\rightarrow\sigma))
\rightarrow((\varphi\rightarrow\psi)\rightarrow(\varphi\rightarrow\sigma))$

(P3) $(\neg\varphi\rightarrow\neg\psi)\rightarrow ((\neg\varphi\rightarrow\psi)\rightarrow\varphi)$,

\vskip 0.1in

\noindent together with the inference rule Modus Ponens ($MP$). So for every $K$, $\{{\rm P1}, {\rm P2}, {\rm P3}\}\subset \textsf{Ax}(K)$ and $MP\in \textsf{IR}(K)$.
In addition each $K$  contains   axioms for the new connective $|$. These  are some  or all of the following schemes.

\vskip 0.1in

($S_1$) \ \ \ $\varphi\wedge \psi\rightarrow \varphi|\psi$

($S_2$) \ \ \ $\varphi|\psi\rightarrow \varphi\vee\psi$

($S_3$) \ \ \ $\varphi|\psi\rightarrow \psi|\varphi$

($S_4$) \ \ \ $(\varphi|\psi)|\sigma\rightarrow \varphi|(\psi|\sigma)$

($S_5$) \ \ \ $\varphi\wedge\neg\psi\rightarrow (\varphi|\psi\leftrightarrow\neg\varphi|\neg\psi)$

\vskip 0.1in

\noindent  Provability (\`{a} la Hilbert) in  $K$, denoted $\vdash_{K}\varphi$,  is defined as usual.  It is clear that
$$\Sigma\vdash \alpha \  \Rightarrow \ \Sigma\vdash_K \alpha,$$
where $\vdash$ denotes  provability in PL.
$\Sigma$ is said to be {\em $K$-consistent,}  if $\Sigma\not\vdash_K\bot$.

Let $K_0$ denote the  formal system described as follows.
$$\textsf{Ax}(K_0)=\{{\rm P1},{\rm P2},{\rm P3}\} +\{S_1,S_2,S_3\},  \quad   \quad \textsf{IR}(K_0)=\{MP\}.$$
Extensions of $K_0$ defined below  will contain also the rule $SV$ (from {\em salva veritate}) defined as follows.
$$(SV) \quad \quad \mbox{\em from} \ \ \varphi\leftrightarrow\psi \ \ \mbox{\em infer} \ \varphi|\sigma\leftrightarrow\psi|\sigma,$$
$$\quad \quad \quad \quad \quad \mbox{if} \ \varphi\leftrightarrow\psi \ \mbox{is provable in $K_0$}.$$
The rule $SV$  guarantees  that if $\alpha$, $\beta$ are classical logically equivalent sentences, then truth is preserved if $\alpha$ is substituted for  $\beta$  in expressions containing $|$ (just as in the case with the standard connectives). Let the formal systems $K_1$, $K_2$ and $K_3$ be defined as follows.
$$\textsf{Ax}(K_1)=\textsf{Ax}(K_0),   \quad \ \quad \textsf{IR}(K_1)=\{\textit{MP},SV\},$$
$$\textsf{Ax}(K_2)=\textsf{Ax}(K_1)+S_4,   \quad \ \quad \textsf{IR}(K_2)=\{\textit{MP},SV\},$$
$$\textsf{Ax}(K_3)=\textsf{Ax}(K_2)+S_5,   \quad \   \quad \textsf{IR}(K_3)=\{\textit{MP},SV\}.$$
A consequence of  $SV$ is that if $\vdash_{K_0}(\varphi\leftrightarrow\psi)$ then, for any $\sigma$, $\vdash_{K_i}(\varphi|\sigma\leftrightarrow\psi|\sigma)$, for $i=1,2,3$.

So much for the syntax of PLS. We now turn to the semantics. The axioms $S_i$  are  motivated by the  intended meaning of $|$ already mentioned above, and the corresponding semantics for sentences of $L_s$ based on choice functions. This semantics consists  of pairs $\langle v,f\rangle$, where $v:Sen(L)\rightarrow \{0,1\}$ is a usual two-valued assignment  of the sentences of $L$, and $f$ is a choice function for pairs of elements of $Sen(L)$, i.e., $f:[Sen(L)]^2\rightarrow Sen(L)$ such that $f(\{\alpha,\beta\})\in \{\alpha,\beta\}$, where for any set $A$, $[A]^2=\{\{a,b\}:a,b\in A\}$. (For basic facts about choice functions the reader may consult \cite{Je73}.) The functions  $f$ are  defined also for singletons with $f(\{\alpha\})=\alpha$. We simplify notation by writing $f(\alpha,\beta)$ instead of $f(\{\alpha,\beta\})$, thus by convention $f(\alpha,\beta)=f(\beta,\alpha)$ and $f(\alpha,\alpha)=\alpha$. $f$ gives rise to a function $\overline{f}:Sen(L_s)\rightarrow Sen(L)$, defined inductively as follows.

\begin{Def} \label{D:collapse}
{\em (i) $\overline{f}(\alpha)=\alpha$, for $\alpha\in Sen(L)$,

(ii) $\overline{f}(\varphi\wedge\psi)=\overline{f}(\varphi)\wedge\overline{f}(\psi)$,

(iii) $\overline{f}(\neg\varphi)=\neg\overline{f}(\varphi)$,

(iv) $\overline{f}(\varphi|\psi)=f(\overline{f}(\varphi),\overline{f}(\psi))$.}
\end{Def}

\noindent We refer to $\overline{f}$ as the {\em collapsing function} induced by $f$. Then we define the truth of $\varphi$ in $\langle v,f\rangle$, denoted $\langle v,f\rangle\models_s\varphi$, as follows.
\begin{equation} \label{E:struth}
\langle v,f\rangle\models_s\varphi: \Leftrightarrow v(\overline{f}(\varphi))=1.
\end{equation}
(In \cite{Tz17} we denote by $M$ the two-valued assignments of sentences of $L$ and write $\langle M,f\rangle$ instead of $\langle v,f\rangle$. Also we write  $M\models \alpha$ instead of  $M(\alpha)=1$.)

We shall refer to the semantics defined by (\ref{E:struth}) as {\em sentence choice semantics,} or SCS for short. A remarkably similar notion of choice function for pairs of sentences, and its interpretation as a ``conservative'' binary connective, was given also independently in \cite{Hum11} (see Example 3.24.14, p. 479).

The reason that we used four formal systems $K_0$-$K_3$,  in increasing strength,  is that they correspond to four different classes of choice functions defined below.

\begin{Def} \label{D:tautologies}
{\em Let ${\cal F}$ denote the set of all choice functions for $Sen(L)$ and let  $X\subseteq {\cal F}$.

(i) For a set $\Sigma\subseteq Sen(L_s)$ and $X\subseteq {\cal F}$, $\Sigma$ is said to be} $X$-satisfiable {\em if there are $v$  and $f\in X$ such that $\langle v,f\rangle \models_s\Sigma$.

(ii) For $\Sigma\subseteq Sen(L_s)$ and $\varphi\in Sen(L_s)$, $\varphi$ is an} $X$-logical consequence of $\Sigma$, {\em denoted $\Sigma\models_X\varphi$, if for every $v$ and every $f\in X$, $\langle v,f\rangle\models_s\Sigma \Rightarrow \langle v,f\rangle\models_s\varphi$.

(iii) $\varphi$ is an} $X$-tautology, {\em denoted $\models_X\varphi$, if $\emptyset\models_X\varphi$.

iv) $\varphi$ and $\psi$ are} $X$-logically equivalent, {\em  denoted $\varphi\sim_X\psi$, if $\models_X(\varphi\leftrightarrow \psi)$. Also let
$$Taut(X)=\{\varphi\in Sen(L_s):\models_X\varphi\}.$$ }
\end{Def}
One of the motivating results behind the development of PLS was the following ``interpolation property'' of $\varphi|\psi$ with respect to $\varphi\wedge\psi$ and $\varphi\vee\psi$ (see Theorem 2.8 of \cite{Tz17}).

\begin{Fac} \label{F:interpol}
For all $\varphi,\psi\in Sen(L_s)$,
$$\varphi\wedge \psi \models_{\cal F}  \varphi|\psi \models_{\cal F}  \varphi\vee \psi,$$
while in general
$$\varphi\vee \psi \not\models_{\cal F}\varphi|\psi \not\models_{\cal F} \varphi\wedge \psi.$$
\end{Fac}

Now while the axioms of $K_0$ are easily seen to be $\models_{\cal F}$-tautologies,  this is not the case  with the  axioms $S_4$ and $S_5$. They correspond to some special subclasses of ${\cal F}$ described below.

\begin{Def} \label{D:stringent}
{\em 1) An  $f\in {\cal F}$ is said to be } associative {\em if for all $\alpha,\beta,\gamma\in Sen(L)$
$$f(f(\alpha,\beta),\gamma)=f(\alpha,f(\beta,\gamma)).$$

2) An $f\in {\cal F}$ is said to be} regular {\em if for all $\alpha,\alpha',\beta\in Sen(L)$,
$$\alpha\sim\alpha' \Rightarrow f(\alpha,\beta)\sim f(\alpha',\beta)$$
where $\alpha\sim\beta$ denotes logical equivalence in PL.
}
\end{Def}
Let
$$\textit{Asso}=\{f\in {\cal F}: f \ \mbox{is asociative}\},$$
$$\textit{Reg}=\{f\in {\cal F}: f \ \mbox{is regular}\},$$

We have the following simple and nice characterization of the functions in $\textit{Asso}$.

\begin{Lem} \label{L:ass}
({\rm \cite[Corollary 2.17]{Tz17})} $f\in Asso$ if and only if there is  a total $<$ ordering of $Sen(L)$ such that $f=\min_<$, i.e., $f(\alpha,\beta)=\min(\alpha,\beta)$ for all $\alpha,\beta\in Sen(L)$.
\end{Lem}
(Actually \ref{L:ass} holds for associative choice functions on an arbitrary set $A$, see Theorem 2.14 of \cite{Tz17}.)
Both properties of associativity and regularity are strongly desirable and would be combined. Also, in view of the above characterization of associative functions through total orderings,  the following definition is natural.

\begin{Def} \label{D:regasso}
{\em A total ordering $<$ of $Sen(L)$ is} regular {\em if the corresponding choice function $f=\min_<$ is regular or, equivalently, if for all $\alpha$, $\beta$ in $Sen(L)$
$$\alpha\not\sim\beta \ \& \ \alpha<\beta \ \Rightarrow \ [\alpha]<[\beta],$$
where $[\alpha]$ is the $\sim$-equivalence class of $\alpha$.
}
\end{Def}
Let
$$\textit{Reg}^*=\textit{Reg}\cap \textit{Asso}.$$
Clearly $f\in \textit{Reg}^*$ iff $f=\min_<$ for a regular total ordering $<$ of $Sen(L)$.

\begin{Def} \label{D:dec}
{\em Let $<$ be a total ordering of $Sen(L)$. $<$ is said to be}  $\neg$-decreasing {\em if for all $\alpha,\beta\in Sen(L)$ such that $\alpha\not\sim\beta$,
$$\alpha<\beta \Leftrightarrow \neg\beta<\neg\alpha.$$
If $f\in \textit{Reg}^*$,  $f$ is said to be} $\neg$-decreasing {\em if $f=\min_<$ for some $\neg$-decreasing $<$.}
\end{Def}
Let
$$\textit{Dec}=\{f\in \textit{Reg}^*: f \ \mbox{is $\neg$-decreasing}\}.$$
Since $\textit{Dec} \subseteq\textit{Reg}^*\subseteq \textit{Reg} \subseteq{\cal F}$, it follows that
$$Taut({\cal F})\subseteq  Taut(Reg)\subseteq Taut(Reg^*)\subseteq Taut(Dec).$$

We can give now a full specification  of the meaning of the notation
$${\rm PLS}(X,K)$$
already introduced in the beginning of this section: given a set $X\subseteq {\cal F}$, and a formal system $K$ with  $\textsf{Ax}(K)\subseteq Taut(X)$, ${\rm PLS}(X,K)$ is the logic with logical  consequence relation $\models_X$, determined by the structures $\langle v,f\rangle$, with $f\in X$, and with provability relation $\vdash_K$.  Given a logic ${\rm PLS}(X,K)$,  the soundness and completeness theorems for it refer as usual to the  connections between the relations $\models_X$ and $\vdash_K$, or between $X$-satisfiability and $K$-consistency.

At this point a word of caution is needed. As is well-known the soundness theorem (ST) and completeness theorem (CT) of a logic have two distinct formulations which are equivalent for classical logic, but need not be so in general.  For the logic ${\rm PLS}(X,K)$ these two forms, ST1 and ST2 for Soundness and CT1 and  CT2 for Completeness, are the following.

$$({\rm ST1}) \hspace{.5\columnwidth minus .5\columnwidth} \Sigma\vdash_K\varphi \ \Rightarrow \ \Sigma\models_X\varphi, \hspace{.5\columnwidth minus .5\columnwidth} \llap{}$$
$$({\rm ST2}) \hspace{.5\columnwidth minus .5\columnwidth}
\Sigma \ \mbox{is $X$-satisfiable} \ \Rightarrow \ \Sigma \ \mbox{is $K$-consistent}
\hspace{.5\columnwidth minus .5\columnwidth} \llap{}$$
$$({\rm CT1}) \hspace{.5\columnwidth minus .5\columnwidth} \Sigma\models_X\varphi \ \Rightarrow \ \Sigma\vdash_K\varphi,  \hspace{.5\columnwidth minus .5\columnwidth} \llap{}$$
$$({\rm CT2}) \hspace{.5\columnwidth minus .5\columnwidth}
\Sigma \ \mbox{is $K$-consistent} \ \Rightarrow \ \Sigma \ \mbox{is $X$-satisfiable}.
\hspace{.5\columnwidth minus .5\columnwidth} \llap{}$$
ST1 and ST2 are easily shown to be  equivalent for every system ${\rm PLS}(X,K)$. Moreover the Soundness Theorem for each one of  the logics ${\rm PLS}({\cal F},K_0)$, ${\rm PLS}(Reg,K_1)$, ${\rm PLS}(Reg^*,K_2)$ and ${\rm PLS}(Dec,K_3)$ is easily established. But
the  equivalence of CT1 and CT2 is based on  the {\em Deduction Theorem} (DT) which is not known to be true for every ${\rm PLS}(X,K)$, when $K$ contains  the inference rule $SV$. Recall that DT is the following implication.  For all $\Sigma$, $\varphi$, $\psi$,
\begin{equation}
\Sigma\cup\{\varphi\}\vdash_K \psi \ \Rightarrow \ \Sigma\vdash_K\varphi\rightarrow\psi.
\end{equation}
Concerning the relationship between  CT1 and CT2 for ${\rm PLS}(X,K)$ the following holds.

\begin{Fac} \label{F:eqsat}
${\rm CT1}\Rightarrow {\rm CT2}$ holds for every ${\rm PLS}(X,K)$. If \ $\vdash_K$ satisfies DT, then the converse holds too, i.e., ${\rm CT1}\Leftrightarrow {\rm CT2}$.
\end{Fac}
The system ${\rm PLS}({\cal F},K_0)$, whose only inference rule is $MP$, satisfies   ${\rm CT1}\Leftrightarrow {\rm CT2}$  as a consequence of DT.  So  we can  just say it is ``complete'' instead of ``CT1-complete'' and ``CT2-complete''. The following is shown in  \cite[\S 3.1]{Tz17})

\begin{Thm} \label{T:mainK0}
${\rm PLS}({\cal F},K_0)$ is complete.
\end{Thm}

However in the  systems over $K_i$, for $i>0$, that contain the extra rule $SV$,  the status of DT is open, so the distinction between CT1 and CT2 remains.
So concerning  the logics ${\rm PLS}(Reg,K_1)$, ${\rm PLS}(Reg^*,K_2)$ and ${\rm PLS}(Dec,K_3)$ it is reasonable to try to prove the weaker of the two forms of completeness, namely  CT2-completeness.  But even this will be proved only conditionally. Because there is still another serious impact  of the lack of DT. This  is  that we don't know if  every consistent set of sentences can be extended to a consistent and {\em complete} set (i.e., one that contains one of the $\varphi$ and  $\neg\varphi$, for every $\varphi$). Of course every consistent set $\Sigma$ can be extended (e.g. by Zorn's Lemma) to a {\em maximal} consistent set $\Sigma'\supseteq\Sigma$. But maximality of $\Sigma'$  does not guarantee completeness without DT. Because  $\Sigma'$ may be maximal consistent and yet there is a $\varphi$ such that $\varphi\notin \Sigma'$ and $\neg\varphi\notin\Sigma'$, so  $\Sigma\cup\{\varphi\}$ and $\Sigma\cup\{\neg\varphi\}$ are both inconsistent. That  looks strange but we don't see how it could be proved false without DT. This property of extendibility of a consistent set to a consistent and complete one, for a formal system $K$, is crucial  for the proof of completeness of $K$ (with respect to a given semantics), so we isolate it as a property of $K$ denoted $cext(K)$. It reads as follows.
$$(cext(K))
\hspace{.5\columnwidth minus .5\columnwidth} \mbox{\em Every $K$-consistent set of sentences can be extended to }\hspace{.5\columnwidth minus .5\columnwidth} \llap{}$$
\hspace{0.9in} $\mbox{\em a $K$-consistent and complete set}$.

\vskip 0.1in

Then the following conditional ${\rm CT2}$-completeness results are shown in  \cite[\S 3.2]{Tz17}).

\begin{Thm} \label{T:main15}

(i) ${\rm PLS}(Reg,K_1)$ is ${\rm CT2}$-complete if and only if  $cext(K_1)$ is true.

(ii) ${\rm PLS}(Reg^*,K_2)$ is ${\rm CT2}$-complete if and only if  $cext(K_2)$ is true.

(iii) ${\rm PLS}(Dec,K_3)$ is ${\rm CT2}$-complete if and only if  $cext(K_3)$ is true.
\end{Thm}

\section{First-order superposition logic and their semantics}

\subsection{What is first-order superposition logic}
First  let  us make precise  what first-order superposition logic, or FOLS for short, is. At axiomatic level the formal systems of FOLS extend the  formal system of first-order logic (FOL) exactly as the formal systems of PLS outlined in section 1.1. extend the  formal system of propositional  logic (PL). So we first fix an axiomatization of  FOL (see e.g. \cite{En02}). If   $L$ is a first-order language with logical constants $\wedge$, $\neg$, $\forall$ and equality $=$, and variables $v_i$, $v$, $u$ etc, let  $L_s=L\cup\{|\}$, where $|$ is the new binary connective for superposition. Let $Fml(L_s)$, $Sen(L_s)$ denote the sets of {\em all} formulas and sentences of $L_s$, respectively,  defined by the usual recursion, as those of $L$, plus the step for the connective $|$.  We stress the word ``all'' because in some  version of formalization considered below, {\em restrictions} to the formation of formulas of $L_s$, concerning the applicability of $|$, might be  sensible.
For example in one of the semantics of FOLS considered below formulas of the form, e.g. $\forall v\exists u(\alpha(v)|\beta(v)|\gamma (u))$ are allowed, while $(\forall v(\alpha(v)|\beta(v)))|(\exists u\gamma (u))$ are not. Thus we reserve the right to deal later only with some subsets of $Fml(L_s)$, $Sen(L_s)$.

The axioms and rules of inference of FOL are the following.
$$\textsf{Ax}({\rm FOL})=\textsf{Ax}({\rm PL})+\{\textit{UI},\textit{D}\}+\{I_1,\ldots,I_5\}, \quad  \textsf{IR}({\rm FOL})=\{\textit{MP}, \textit{GR}\},$$
where $GR$ is the generalization rule, $UI$ (Universal Instantiation scheme) and $D$ are the basic axioms of FOL (for the language $L_s$) concerning quantifiers, $I_1$, $I_2$, $I_3$ are the trivial axioms for $=$ (reflection, symmetry and transitivity),  and finally $I_4$ and $I_5$ are  the schemes of substitution of equals within  terms and formulas. Specifically:

({\textit{UI}) \ \ \ $\forall v\varphi(v) \rightarrow \varphi(t)$, for every closed term $t$,

($D$) \ \ \ $\forall v(\varphi\rightarrow\psi(v))\rightarrow (\varphi\rightarrow \forall v\psi(v))$, if $v$ is not free in $\varphi$,

($I_4$) \ \ \ $(\forall v,u)(v=u \rightarrow t(v)=t(u))$,

($I_5$) \ \ \  $(\forall v,u)(v=u \rightarrow (\varphi(v)\rightarrow \varphi(u))$.

\noindent   The sentences of $L$ will be interpreted in $L$-structures ${\cal M}=\langle M,\ldots\rangle$.

\vskip 0.1in

{\bf Notational convention.} We  keep using the notational convention introduced in the previous section,  that  throughout the letters $\varphi$, $\psi$, $\sigma$ will denote in general formulas  of  $L_s$, while the letters $\alpha$, $\beta$, $\gamma$ are reserved for  formulas of  $L$ only.

The extra axioms of FOLS will be among the schemes already seen in 1.1, namely:

\vskip 0.1in

($S_1$) \ \ \ $\varphi\wedge \psi\rightarrow \varphi|\psi$

($S_2$) \ \ \ $\varphi|\psi\rightarrow \varphi\vee\psi$

($S_3$) \ \ \ $\varphi|\psi\rightarrow \psi|\varphi$

($S_4$) \ \ \ $(\varphi|\psi)|\sigma\rightarrow \varphi|(\psi|\sigma)$

($S_5$) \ \ \ $\varphi\wedge\neg\psi\rightarrow (\varphi|\psi\leftrightarrow\neg\varphi|\neg\psi)$

\vskip 0.1in

The formal systems we are going to deal with below are $\Lambda_0$, $\Lambda_1$,  $\Lambda_2$ and  $\Lambda_3$ defined as follows

$$\textsf{Ax}(\Lambda_0)=\textsf{Ax}({\rm FOL})+\{S_1,S_2,S_3\}, \quad \textsf{IR}(\Lambda_0)=\{MP,GR\}$$
$$\textsf{Ax}(\Lambda_1)=\textsf{Ax}(\Lambda_0), \quad \textsf{IR}(\Lambda_1)=\{MP,GR,SV\}$$
$$\textsf{Ax}(\Lambda_2)=\textsf{Ax}(\Lambda_1)+S_4, \quad \textsf{IR}(\Lambda_2)=\{MP,GR, SV\},$$
$$\textsf{Ax}(\Lambda_3)=\textsf{Ax}(\Lambda_2)+S_5, \quad \textsf{IR}(\Lambda_3)=\{MP,GR, SV\},$$
where $SV$ is the rule Salva Veritate mentioned in section 1.1, but with $\Lambda_0$ in place of $K_0$. That is:
$$(SV) \quad \quad \mbox{\em from} \ \ \varphi\leftrightarrow\psi \ \ \mbox{\em infer} \ \varphi|\sigma\leftrightarrow\psi|\sigma,$$
$$\quad \quad \quad \quad \quad \mbox{if} \ \varphi\leftrightarrow\psi \ \mbox{is provable in $\Lambda_0$}.$$
Note that there are no new axioms for $|$ in FOLS beyond those of PLS, which means that there is no natural interplay between $|$ and quantifiers. In fact the  connections one might consider between $|$ and $\forall$, e.g. $(\forall v)(\varphi|\psi)\leftrightarrow (\forall v \varphi)|(\forall v\psi)$, or $(\exists v)(\varphi|\psi)\leftrightarrow (\exists v \varphi)|(\exists v\psi)$,  either do not make sense because the formulas involved are illegitimate, or are simply false in the  semantics where the formulas involved are allowed (e.g. in the semantics of \cite{Tz18}).

\subsection{Candidate semantics for FOLS}
In \cite{Tz18} we developed an alternative semantics (and a slightly different formalization) for PLS, based on choice function not for pairs of sentences but for pairs of elements of a Boolean algebra ${\cal B}$ where the classical sentences take truth values. We called this ``Boolean-value choice semantics'', or BCS for short. It turned out that this semantics can apply also to FOLS without extra pains, and with respect to this semantics the formal systems of FOLS satisfy some natural soundness and completeness results.

The main question addressed in this paper is  whether FOLS can admit a semantics that naturally  extends and generalizes  the sentence choice semantics (SCS) of \cite{Tz17} (based on the truth definition (\ref{E:struth}) mentioned in the previous section). ``Naturally'' means that the semantics will continue to consist of  pairs $\langle {\cal M},f\rangle$, where ${\cal M}$ is an  $L$-structure and $f$ is  a choice function for pairs of formulas/sentences of $L$, and will  follow  the basic reduction of truth to the Tarskian one  through the relation:  $\langle {\cal M},f\rangle\models_s\varphi $ iff ${\cal M}\models \overline{f}(\varphi)$. The intricate question is about the  {\em domain} of the choice function $f$.  Namely, would  $f$  apply to  pairs of {\em all} formulas of $L$, or only to {\em some} such pairs, e.g. to pairs of sentences alone? The answer to the above question is that there can be no natural extension of SCS to a ``formula choice semantics'', in the sense that $f$ is allowed to apply to pairs of {\em arbitrary}  formulas. Such a semantics fails badly for reasons independent of the connective $|$, simply as a result of incompatibility between  choice of formulas with free variables and corresponding choice of formulas with substituted terms. On the other hand, it is shown that  a semantics with some restrictions both to the construction of formulas, as well as to the applicability of choice functions (allowing them to apply to pairs of sentences only), can work smoothly and lead to satisfactory soundness and completeness results with respect to the axiomatization of FOLS, which is essentially the same as the one considered in \cite{Tz18}.

First let us note that in  any case, whatever the domain of $f$ would be, the collapsing map $\overline{f}$ should satisfy conditions (i)-(iv) of Definition  \ref{D:collapse}. So if we assume that $f$ is defined for all pairs of  quantifier-free sentences of  $L$,  then  conditions (i)-(iv), in combination with the truth definition (\ref{E:struth}), suffice  to define $\langle {\cal M},f\rangle\models_s\varphi$ for every quantifier-free sentence of  $L_s$. So the only missing  step for the complete definition of $\langle {\cal M},f\rangle\models_s\varphi$ is the definition of  $\langle {\cal M},f\rangle\models_s\forall v\varphi(v)$. For that we have two options, called  {\em formula choice semantics} (FCS for short) and  {\em sentence choice semantics} (SCS) because they are based on the use of choice functions for pairs of arbitrary formulas and for pairs of sentences alone, respectively. To distinguish them we shall use the  symbols $\models_s^1$ and $\models_s^2$ for the resulting  truth relations, respectively.

\vskip 0.1in

{\bf Option 1. Formula choice semantics (FCS)} Here the set of formulas $Fml(L_s)$ of $L_s$ is defined by the usual closure steps with respect to the connectives (including $|$) and quantifiers and every choice function $f$ is defined on the entire $[Fml(L)]^2$. Therefore the truth definition of quantified sentences  should be as follows.
\begin{equation} \label{E:basictruth}
\langle {\cal M},f\rangle\models^1_s\forall v\varphi(v)\Leftrightarrow  {\cal M}\models\overline{f}(\forall v\varphi(v)).
\end{equation}
It is easy to see that this definition is meaningful and effective if and only if the collapsing mapping $\overline{f}$ commutes with $\forall$, i.e., if  $\overline{f}$ satisfies, in addition to conditions (i)-(iv) of Definition \ref{D:collapse}, the condition:

\vskip 0.1in

(v) $\overline{f}(\forall v\varphi)=\forall v \overline{f}(\varphi)$.\footnote{Otherwise, one cannot see how e.g. $\overline{f}(\forall v(\alpha|\beta))$ could be defined.}

\vskip 0.1in

\noindent  [Treating $\exists$ as usual, i.e., as $\neg\forall\neg$, it  follows from (v) and (iii) of \ref{D:collapse} that $\overline{f}(\exists v\varphi)=\exists v \overline{f}(\varphi)$.]

Throughout this subsection we shall often refer to conditions (i)-(iv) of \ref{D:collapse} together with condition (v) above as ``conditions (i)-(v)'' for $\overline{f}$.   By (v),  (\ref{E:basictruth}) becomes
\begin{equation} \label{E:option1}
\langle {\cal M},f\rangle\models^1_s\forall v\varphi(v)\Leftrightarrow  {\cal M}\models\forall v\overline{f}(\varphi(v)).
\end{equation}
The right-hand side of (\ref{E:option1}) is an instance of  Tarskian satisfaction, so it holds iff ${\cal M}\models\overline{f}(\varphi(v))(x)$ is true for every $x\in M$, where the elements of $M$ are used as parameters added to $L$. Therefore (\ref{E:option1}) is equivalently written
\begin{equation} \label{E:option10}
\langle {\cal M},f\rangle\models^1_s\forall v\varphi(v)\Leftrightarrow  {\cal M}\models\overline{f}(\varphi(v))(x), \ \mbox{for every $x\in M$}.
\end{equation}
Thus  (\ref{E:option1}) (or (\ref{E:option10}))  determines the truth of every sentence $\varphi$ of $L_s(M)$ in $\langle {\cal M},f\rangle$ with respect to $\models_s^1$. We refer to the truth relation $\models_s^1$ (for obvious reasons) as {\em formula choice semantics}, or FCS for short. We shall see however below that  FCS fails badly not with respect to the interpretation of $|$, but because, surprisingly enough, fails to satisfy the Universal Instantiation scheme, as a consequence of the fact that $f$ applies to pairs of formulas with free variables. So a reasonable alternative would be to restrict $f$ to pairs of sentences alone.

\vskip 0.1in

{\bf Option 2. Sentence choice semantics (SCS)} Assume now that the choice functions $f$ are defined only for {\em sentences} of $L(M)=L\cup M$, where the latter is $L$ augmented with the elements of $M$ treated as parameters. We let the letters $x,y,a,c$ range over elements of $M$. The question is how the collapsing  $\overline{f}$ is defined in this case and for which $\varphi$ of $L_s$. For instance, what would  $\overline{f}(\forall v(\alpha(v)|\beta(v)))$ be for classical $\alpha(v)$ and $\beta(v)$? Letting  $\overline{f}(\forall v(\alpha(v)|\beta(v)))=\forall v\overline{f}(\alpha(v)|\beta(v))=\forall vf(\alpha(v),\beta(v))$ is not an option since $f$ does not apply to pairs of open formulas. The answer is simply that for $Q\in\{\forall,\exists\}$,
\begin{equation} \label{E:notdef}
\overline{f}(Qv(\alpha(v)|\beta(v))) \  \mbox{are not defined}.
\end{equation}
However this is not necessarily a dead end. It would only prompt us to define the truth of $\forall v(\alpha(v)|\beta(v))$ in $\langle {\cal M},f\rangle$ not through  (\ref{E:option1}), but in the  Tarskian way:
$$\langle {\cal M},f\rangle\models_s\forall v(\alpha(v)|\beta(v)))\Leftrightarrow \langle {\cal M},f\rangle\models_s(\alpha(x)|\beta(x)),$$
for all $x\in M$. So let us define, alternatively to (\ref{E:option1}),  for every universal well-formed formula  $\forall v\varphi(v)$ of $L_s$:
\begin{equation} \label{E:option21}
\langle {\cal M},f\rangle\models^2_s\forall v\varphi(v)\Leftrightarrow \langle{\cal M},f\rangle\models_s^2\varphi(x), \ \mbox{for every $x\in M$}.
\end{equation}
From (\ref{E:option21}), combined with clause (iii) of \ref{D:collapse}, clearly we have  also that
\begin{equation} \label{E:option22}
\langle {\cal M},f\rangle\models^2_s\exists v\varphi(v)\Leftrightarrow \langle{\cal M},f\rangle\models_s^2\varphi(x), \ \mbox{for some $x\in M$}.
\end{equation}
(\ref{E:option21}) and (\ref{E:option22}) settle the definition with respect to $\models_s^2$ of sentences that begin with a quantifier. This also implicitly suggests that for $\varphi$ that do not begin with a quantifier, the truth of $\varphi$ in $\langle {\cal M},f\rangle$ should be  defined by means of the collapsing map $\overline{f}$ i.e.,
\begin{equation} \label{E:option23}
\langle {\cal M},f\rangle\models^2_s\varphi\Leftrightarrow {\cal M}\models\overline{f}(\varphi).
\end{equation}
But this will immediately lead to trouble, unless we  put restrictions to the formation of formulas of $L_s$. For consider, say,  the sentence  $(\forall v(\alpha|\beta))|(\exists u(\gamma|\delta))$. Then we should have
$$\langle {\cal M},f\rangle\models^2_s(\forall v(\alpha|\beta))|(\exists u(\gamma|\delta))\Leftrightarrow {\cal M}\models\overline{f}[(\forall v(\alpha|\beta))|(\exists u(\gamma|\delta))]\Leftrightarrow$$
$${\cal M}\models f(\overline{f}(\forall v(\alpha|\beta)),\overline{f}(\exists u(\gamma|\delta))).$$
But by (\ref{E:notdef}) above, $\overline{f}(\forall v(\alpha|\beta))$ and $\overline{f}(\exists u(\gamma|\delta))$ are not defined, so the last part of the above equivalences  does not make sense.

The conclusion is that if we want to employ choice functions for pairs of {\em sentences} only and $\models_s^2$ obeys (\ref{E:option21}) and (\ref{E:option22}), formulas like $(\forall v(\alpha|\beta))|(\exists u(\gamma|\delta))$, should not be allowed. That is, instead  of the full set of formulas $Fml(L_s)$ we shall consider the restricted set of formulas $RFml(L_s)$. The latter differs from $Fml(L_s)$ in that  $\varphi|\psi$ belongs to $RFml(L_s)$ iff $\varphi$ and $\psi$ are either classical or quantifier free.  We shall refer to the truth relation $\models_s^2$ as {\em sentence  choice semantics}, or SCS, just as we did with the corresponding semantics of PLS. We shall examine $\models_s^2$ in more detail in section 5.

\vskip 0.1in

Obviously the two semantics based on  $\models_s^1$ and $\models_s^2$ are  not equivalent, since they apply to different sets of sentences. However, there are sentences $\varphi$ for which  both truth definitions $\langle {\cal M},f\rangle\models^1_s\varphi$  and $\langle {\cal M},f\rangle\models^2_s\varphi$ make sense.  But in general even for such sentences the definitions  do not coincide.  Actually none of them implies the other. The difference is easily detected by observing the right-hand sides of (\ref{E:option10}) and (\ref{E:option21}) for $\varphi$ that do not begin with a quantifier. Namely, $\overline{f}(\varphi(v))(x)$ and $\overline{f}(\varphi(x))$ are in general  inequivalent. To illustrate it, let $\varphi(v):=\alpha(v)|\beta(v)$, where $\alpha(v)$ and $\beta(v)$ are formulas of $L$.  Let ${\cal M}$ be an $L$-structure. To compare the two approaches, we must use an $f$ which is meaningful in both of them, i.e., it applies to all pairs of formulas of $L(M)$.  Fix such an $f$. Then $f$ applies to  $\{\alpha(v),\beta(v)\}$ and let $f(\alpha(v),\beta(v))=\alpha(v)$. By (\ref{E:option10}),  $\langle{\cal M},f\rangle\models^1_s\forall v(\alpha(v)|\beta(v))$ iff  ${\cal M}\models  f(\alpha(v),\beta(v))(x)$, for all $x\in M$, therefore
\begin{equation} \label{E:ex1}
\langle{\cal M},f\rangle\models^1_s\forall v(\alpha(v)|\beta(v)) \Leftrightarrow {\cal M}\models\alpha(x), \ \mbox{for all $x\in M$}.
\end{equation}
On the other hand,   by (\ref{E:option21}),  $\langle{\cal M},f\rangle\models^2_s\forall v(\alpha(v)|\beta(v))$ iff  $\langle{\cal M},f\rangle\models^2_s \alpha(x)|\beta(x)$, for all $x\in M$, hence
\begin{equation} \label{E:ex2}
\langle{\cal M},f\rangle\models^2_s\forall v(\alpha(v)|\beta(v))\Leftrightarrow {\cal M}\models f(\alpha(x),\beta(x)), \ \mbox{for all $x\in M$}.
\end{equation}
The right-hand sides of (\ref{E:ex1}) and (\ref{E:ex2}) may be quite different, since the  choices of $f$ from the pairs $\{\alpha(x),\beta(x)\}$, for the various $x\in M$, may be non-uniform, e.g. for $x_1\neq x_2$  we may have $f(\alpha(x_1),\beta(x_1))=\alpha(x_1)$ and $f(\alpha(x_2),\beta(x_2))=\beta(x_2)$. In order for the definitions (\ref{E:ex1}) and (\ref{E:ex2}) to be equivalent, $f$ should be a {\em uniform choice function}, i.e., $f(\alpha(\vec{v}),\beta(\vec{v}))=\alpha(\vec{v})$ should imply  $f(\alpha(\vec{t}),\beta(\vec{t}))=\alpha(\vec{t})$ for all pairs of formulas $\{\alpha(\vec{v}),\beta(\vec{v})\}$ and every tuple of terms $\vec{t}$ that can be substituted for $\vec{v}$. However, as we shall prove in section 4, no choice function $f:[Fml(L)]^2\rightarrow Fml(L)$ can have this property.

\section{The formula choice semantics (FCS) and the failure of universal instantiation}
Since FOLS extends FOL, any proper  semantics for FOLS should first of all satisfy the quantifier axioms of FOL, namely  $\textit{UI}$ and $D$.  In this section we show that unfortunately (and rather unexpectedly)  FCS fails to satisfy $UI$.

Firstly recall that given a language $L$,  whenever we write $\varphi(\vec{v})$, for a formula of $L_s$,  we mean that the free variables of $\varphi$ are {\em among} those of the tuple $\vec{v}$. Then  the following can be  easily verified  by induction on the length of $\varphi$.

\begin{Fac} \label{F:variables}
For every choice function $f$ for pairs of formulas and every $\varphi\in Fml(L_s)$, the free variables of $\overline{f}(\varphi)$ are included in those of $\varphi$, i.e., $FV(\overline{f}(\varphi))$ $\subseteq FV(\varphi)$. In particular, if $\varphi$ is a sentence of $L_s$, then  $\overline{f}(\varphi)$ is a sentence of $L$.
\end{Fac}
[In general, $FV(\overline{f}(\varphi))\varsubsetneq FV(\varphi)$, since, for example, we may have $\varphi(v_1,v_2)=\alpha(v_1)|\beta(v_2)$ and $f(\alpha(v_1),\beta(v_2))=\alpha(v_1)$, so $\overline{f}(\varphi)=f(\alpha(v_1),\beta(v_2))=\alpha(v_1)$.]
It follows from this  Fact that the variables of $\overline{f}(\varphi)$ are  among the variables of $\varphi$, so we may write  for every $\varphi(\vec{v})$:
\begin{equation} \label{E:var}
\overline{f}(\varphi(\vec{v}))=\overline{f}(\varphi)(\vec{v}).
\end{equation}

\begin{Fac} \label{F:axA}
The scheme $D$ is a tautology with respect to FCS.
\end{Fac}

{\em Proof.} Take an instance of $D$
$$\sigma:(\forall v)(\varphi\rightarrow \psi(v))\rightarrow (\varphi\rightarrow (\forall v)\psi(v)),$$  where $\varphi$ does not contain  $v$ free, and take an arbitrary choice function satisfying  conditions  (i)-(v). Then clearly applying these conditions  we have
$$\overline{f}(\sigma)=[(\forall v)(\overline{f}(\varphi)\rightarrow \overline{f}(\psi(v)))\rightarrow (\overline{f}(\varphi)\rightarrow (\forall v)\overline{f}(\psi(v)))].$$
Let  $\overline{f}(\varphi)=\alpha$ and $\overline{f}(\psi(v))=\beta(v)$. Then
the last formula is written
$$\overline{f}(\sigma)=[(\forall v)(\alpha\rightarrow \beta(v))\rightarrow (\alpha\rightarrow (\forall v)\beta(v))].$$
By assumption $v\not\in FV(\varphi)$, and by Fact \ref{F:variables}, $FV(\overline{f}(\varphi))\subseteq FV(\varphi)$, so $v\notin FV(\alpha)$, therefore  $\overline{f}(\sigma)$ is an instance of the scheme $D$ of  FOL, so it holds in every $L$-structure ${\cal M}$. Therefore ${\cal M}\models \overline{f}(\sigma)$, or equivalently $\langle {\cal M}, f\rangle\models \sigma$. \telos

\vskip 0.2in

However the situation is quite different for the scheme $\textit{UI}$. The next theorem shows that, under mild conditions for $L$, there is  no choice function $f$ for $L$ with respect to which $\textit{UI}$  could be a scheme of tautologies.

\begin{Thm} \label{T:nohope}
Let $L$ be a language with at least two distinct closed terms $t_1$, $t_2$, and a formula $\alpha(v)$ in one free variable such that both $\alpha(t_1)\wedge\neg\alpha(t_2)$ and $\neg\alpha(t_1)\wedge\alpha(t_2)$ are satisfiable.  Then for every choice function $f$ for $Fml(L)$ there is an $L$-structure ${\cal M}$ and a formula $\psi(v_1,v_2)$ of $L_s$ with two free variables for which $\textit{UI}$ fails in $\langle {\cal M},f\rangle$, i.e., such that $\langle {\cal M},f\rangle\models^1_s (\forall v_1,v_2)\psi(v_1,v_2)\wedge \neg\psi(t_1,t_2)$.
\end{Thm}

{\em Proof.} [Note first that the conditions required for $L$ in the above Lemma are quite weak. E.g. any $L$ containing three distinct constants $c_1,c_2,c_3$ satisfies them. For if we set $\alpha(v)=(v=c_3)$, then $\alpha(c_1)\wedge\neg\alpha(c_2)$ and $\neg\alpha(c_1)\wedge\alpha(c_2)$ are both satisfiable in $L$-structures.]

Now let $L$, $t_1$, $t_2$ and $\alpha(v)$ be  as stated. Then there are structures ${\cal M}_1$, ${\cal M}_2$ such that ${\cal M}_1\models\alpha(t_1)\wedge\neg\alpha(t_2)$ and  ${\cal M}_2\models\neg\alpha(t_1)\wedge\alpha(t_2)$. Pick a choice function $f$ for $Fml(L)$. It suffices to show that there is a formula $\psi(v_1,v_2)$ of $L_s$   such that either $\langle {\cal M}_1,f\rangle\models^1_s (\forall v_1,v_2)\psi(v_1,v_2)\wedge \neg\psi(t_1,t_2)$, or $\langle {\cal M}_2,f\rangle\models^1_s (\forall v_1,v_2)\psi(v_1,v_2)\wedge \neg\psi(t_1,t_2)$. Let $\alpha(v_2)$ be the formula resulting from $\alpha(v_1)$ if we replace $v_1$ by the new variable $v_2$. We examine how $f$ acts on the pairs of formulas $\{\alpha(v_1),\alpha(v_2)\}$ and $\{\alpha(t_1),\alpha(t_2)\}$ and consider the four possible  cases.

\vskip 0.1in

{\em Case 1.}  $f(\alpha(v_1),\alpha(v_2))=\alpha(v_1)$ and  $f(\alpha(t_1),\alpha(t_2))=\alpha(t_1)$.

By assumption ${\cal M}_2\models\neg\alpha(t_1)\wedge\alpha(t_2)$. Arguing in the standard logic FOL, this can be  written as follows:
$${\cal M}_2\models(\forall v_1 v_2)[v_1=t_2 \wedge v_2=t_1\rightarrow \alpha(v_1)]\wedge \neg[t_2=t_2\wedge t_1=t_1\rightarrow \alpha(t_1)],$$
or
$${\cal M}_2\models(\forall v_1 v_2)[v_1=t_2 \wedge v_2=t_1\rightarrow f(\alpha(v_1),\alpha(v_2))]\wedge$$ $$\neg[t_2=t_2\wedge t_1=t_1\rightarrow f(\alpha(t_1),\alpha(t_2))],$$
or
$$\langle{\cal M}_2,f\rangle\models^1_s(\forall v_1 v_2)[v_1=t_2 \wedge v_2=t_1\rightarrow \alpha(v_1)|\alpha(v_2)]\wedge$$ $$\neg[t_2=t_2\wedge t_1=t_1\rightarrow \alpha(t_1)|\alpha(t_2)],$$
or, since $|$ is commutative,
$$\langle{\cal M}_2,f\rangle\models^1_s(\forall v_1 v_2)[v_1=t_2 \wedge v_2=t_1\rightarrow \alpha(v_1)|\alpha(v_2)]\wedge$$ $$\neg[t_2=t_2\wedge t_1=t_1\rightarrow \alpha(t_2)|\alpha(t_1)].$$
Setting
$$\psi(v_1,v_2):=[v_1=t_2 \wedge v_2=t_1\rightarrow \alpha(v_1)|\alpha(v_2)],$$
the last relation is written
$$\langle{\cal M}_2,f\rangle\models^1_s(\forall v_1 v_2)\psi(v_1,v_2) \wedge\neg\psi(t_2,t_1),$$
thus $\textit{UI}$ fails in $\langle{\cal M}_2,f\rangle$.

\vskip 0.1in

{\em Case 2.}  $f(\alpha(v_1),\alpha(v_2))=\alpha(v_1)$ and  $f(\alpha(t_1),\alpha(t_2))=\alpha(t_2)$.

Now we use the fact that ${\cal M}_1\models\alpha(t_1)\wedge\neg\alpha(t_2)$. As before this is written equivalently,
$${\cal M}_1\models(\forall v_1 v_2)[v_1=t_1 \wedge v_2=t_2\rightarrow \alpha(v_1)]\wedge \neg[t_1=t_1\wedge t_2=t_2\rightarrow \alpha(t_2)],$$
or
$${\cal M}_1\models(\forall v_1 v_2)[v_1=t_1 \wedge v_2=t_2\rightarrow f(\alpha(v_1),\alpha(v_2))]\wedge$$ $$\neg[t_1=t_1\wedge t_2=t_2\rightarrow f(\alpha(t_1),\alpha(t_2))],$$
or
$$\langle{\cal M}_1,f\rangle\models^1_s(\forall v_1 v_2)[v_1=t_1 \wedge v_2=t_2\rightarrow \alpha(v_1)|\alpha(v_2)]\wedge$$ $$\neg[t_1=t_1\wedge t_2=t_2\rightarrow \alpha(t_1)|\alpha(t_2)].$$
Thus putting
$$\psi(v_1,v_2):=[v_1=t_1 \wedge v_2=t_2\rightarrow \alpha(v_1)|\alpha(v_2)],$$
we are done.

\vskip 0.1in

{\em Case 3.} $f(\alpha(v_1),\alpha(v_2))=\alpha(v_2)$ and  $f(\alpha(t_1),\alpha(t_2))=\alpha(t_1)$.

We use again the fact that ${\cal M}_2\models\neg\alpha(t_1)\wedge \alpha(t_2)$, which yields as before
$${\cal M}_2\models(\forall v_1 v_2)[v_1=t_1 \wedge v_2=t_2\rightarrow \alpha(v_2)]\wedge \neg[t_1=t_1\wedge t_2=t_2\rightarrow \alpha(t_1)],$$
or
$$\langle{\cal M}_2,f\rangle\models^1_s(\forall v_1 v_2)[v_1=t_1 \wedge v_2=t_2\rightarrow \alpha(v_1)|\alpha(v_2)]\wedge$$ $$\neg(t_1=t_1\wedge t_2=t_2\rightarrow \alpha(t_1)|\alpha(t_2)].$$
So setting
$$\psi(v_1,v_2):=[v_1=t_1 \wedge v_2=t_2\rightarrow \alpha(v_1)|\alpha(v_2)]$$
we are done.

\vskip 0.1in

{\em Case 4.} $f(\alpha(v_1),\alpha(v_2))=\alpha(v_2)$ and  $f(\alpha(t_1),\alpha(t_2))=\alpha(t_2)$.

We use the fact that ${\cal M}_1\models \alpha(t_1)\wedge\neg\alpha(t_2)$ which translates into
$${\cal M}_1\models(\forall v_1 v_2)[v_1=t_2 \wedge v_2=t_1\rightarrow \alpha(v_2)]\wedge \neg[t_2=t_2\wedge t_1=t_1\rightarrow \alpha(t_2)],$$
or
$$\langle{\cal M}_1,f\rangle\models^1_s(\forall v_1 v_2)[v_1=t_2 \wedge v_2=t_1\rightarrow \alpha(v_1)|\alpha(v_2)]\wedge$$ $$\neg(t_2=t_2\wedge t_1=t_1\rightarrow \alpha(t_2)|\alpha(t_1)].$$
Setting
$$\psi(v_1,v_2):=[v_1=t_2 \wedge v_2=t_1\rightarrow \alpha(v_1)|\alpha(v_2)]$$
we are done. This completes the proof. \telos

\vskip 0.2in

Equivalent to $\textit{UI}$ (in FOL, hence also in FOLS) is the dual axiom of {\em Existential Generalization ($\textit{EG}$):}

\vskip 0.1in

$(\textit{EG}) \ \ \  \varphi(t)\rightarrow (\exists v)\varphi(v)$.

\vskip 0.1in

\noindent Therefore Theorem \ref{T:nohope} is equivalently formulated as follows.

\begin{Cor} \label{C:EG}
Given any language $L$ as above, for every choice function $f$ for $Fml(L)$ there is ${\cal M}$ such that $\textit{EG}$ fails in $\langle {\cal M},f\rangle$, namely, there is a formula $\varphi(v_1,v_2)$ and  closed terms $t_1,t_2$ such that $\langle {\cal M},f\rangle\models_s^1 \varphi(t_1,t_2)\wedge \neg(\exists v_1,v_2)\varphi(v_1,v_2)$.
\end{Cor}

In Theorem \ref{T:nohope} the  formula(s) $\psi(v_1,v_2)$  used to refute $\textit{UI}$  contain   two free variables. We do not know if it possible to refute $\textit{UI}$  using a formula with a single free variable.

Also in the  proof of \ref{T:nohope} we used  a superposed formula of the form $\alpha(v_1)|\alpha(v_2)$, which looks somewhat artificial.  Can we show the failure of $UI$, using a superposition of the form $\alpha(\vec{v})|\beta(\vec{v})$ where $\alpha(\vec{v})$ and $\beta(\vec{v})$ are distinct formulas? The answer is yes. Specifically, by essentially the same argument we can prove the following  variant of  Theorem  \ref{T:nohope}.

\begin{Thm} \label{T:gennohope}
Let $L$ be a language and  assume that  there exist formulas $\alpha(\vec{v})$, $\beta(\vec{v})$  and corresponding tuples of closed terms $\vec{t}$, $\vec{s}$ such that:

(a) $\alpha(\vec{s})=\beta(\vec{t})$,

(b) $\alpha(\vec{t})=\beta(\vec{s})$

(c)  $\alpha(\vec{t})\wedge \neg\beta(\vec{t})$ and $\neg\alpha(\vec{t})\wedge \beta(\vec{t})$ are satisfiable.

\noindent Then for every choice function $f$ there is a structure ${\cal M}$,  a formula $\psi(\vec{v})$ and closed terms $\vec{t}$ such that $\langle {\cal M},f\rangle\models^1_s(\forall\vec{v})\psi(\vec{v}) \wedge \neg \psi(\vec{t})$.
\end{Thm}

{\em Proof.}  [\noindent For example if $<$ is a binary relation of $L$, and $\alpha(v_1,v_2):=(v_1<v_2)$, $\beta(v_1,v_2):=(v_1>v_2)$, $\vec{t}=\langle t_1,t_2\rangle$ and $\vec{s}=\langle t_2,t_1\rangle$, then $\alpha$, $\beta$, $\vec{t}$ and $\vec{s}$ satisfy conditions (a)-(c) above.]

The argument goes exactly as in the proof of Theorem \ref{T:nohope}.
We fix structures ${\cal M}_1$, ${\cal M}_2$ such that ${\cal M}_1\models\alpha(\vec{t})\wedge\neg\beta(\vec{t})$ and  ${\cal M}_2\models\neg\alpha(\vec{t})\wedge\beta(\vec{t})$. Pick a choice function $f$ for $Fml(L)$. It suffices to show that there is a formula $\psi(\vec{v})$ of $L_s$   such that either $\langle {\cal M}_1,f\rangle\models^1_s (\forall \vec{v})\psi(\vec{v})\wedge \neg\psi(\vec{r})$, or $\langle {\cal M}_2,f\rangle\models^1_s (\forall \vec{v})\psi(\vec{v})\wedge \neg\psi(\vec{r})$, for $\vec{r}=\vec{t}$ or $\vec{r}=\vec{s}$.   As before we examine how $f$ acts on the pairs of formulas $\{\alpha(\vec{v}),\beta(\vec{v})\}$, and $\{\alpha(\vec{t}),\beta(\vec{t})\}=\{\alpha(\vec{s}),\beta(\vec{s})\}$, and we examine the four possible cases that  arise as before. Namely:

\vskip 0.1in

{\em Case 1.}  $f(\alpha(\vec{v}),\beta(\vec{v}))=\alpha(\vec{v})$ and  $f(\alpha(\vec{t}),\beta(\vec{t}))=\alpha(\vec{t})$.

\vskip 0.1in

{\em Case 2.}  $f(\alpha(\vec{v}),\beta(\vec{v}))=\alpha(\vec{v})$ and  $f(\alpha(\vec{t}),\beta(\vec{t}))=\beta(\vec{t})$.

\vskip 0.1in

{\em Case 3.}  $f(\alpha(\vec{v}),\beta(\vec{v}))=\beta(\vec{v})$ and  $f(\alpha(\vec{t}),\beta(\vec{t}))=\alpha(\vec{t})$.

\vskip 0.1in

{\em Case 4.}  $f(\alpha(\vec{v}),\beta(\vec{v}))=\beta(\vec{v})$ and  $f(\alpha(\vec{t}),\beta(\vec{t}))=\beta(\vec{t})$.

\vskip 0.1in

\noindent In each of these cases we work as in the corresponding case of the proof of  \ref{T:nohope}.  Details are left to the reader.  \telos

\vskip 0.2in

Closing this section, let us remark that the failure of $UI$  is
{\em fatal} for any quantified logical system like $\Lambda$, in the sense that there can be no reasonable ``weakening'' of $\Lambda$  in which $\forall$ is still in use while $UI$ fails. For the failure of $UI$ is  quite different from the failure e.g. of the Excluded Middle (EM), which has led to a logic weaker than the  classical one and yet quite interesting. The reason is that  $UI$ expresses exactly  the  meaning of ``all'', as a  fundamental logical constant, while  EM does not express the meaning of any logical constant.

\section{The impossibility of uniform choice functions}
At the end of section 2.2, comparing the truth relations $\models_s^1$ and $\models_s^2$,  we said  that the two notions of truth deviate even for  choice functions $f$ that are defined in both semantics,  because $f$ cannot be {\em uniform} when considered as a choice function in FCS. Let us make this claim  precise.

\begin{Def} \label{D:uniff}
{\em A choice function $f:[Fml(L)]^2\rightarrow Fml(L)$ is said to be} uniform {\em
if for any two formulas $\alpha(\vec{v})$, $\beta(\vec{v})$, with $\vec{v}$ free, and any tuple $\vec{t}$ of terms  substitutable for $\vec{v}$ in $\alpha,\beta$, $f(\alpha(\vec{v}),\beta(\vec{v}))\sim \alpha(\vec{v})$ implies $f(\alpha(\vec{t}),\beta(\vec{t}))\sim \alpha(\vec{t})$, or equivalently,  if
the following equivalence holds:
\begin{equation} \label{E:crucial}
[f(\alpha(\vec{v}),\beta(\vec{v}))](\vec{t})\sim f(\alpha(\vec{t}),\beta(\vec{t})).
\end{equation}
}
\end{Def}

{\bf Note.} The reason for writing $\sim$ instead of $=$ in condition  (\ref{E:crucial}) above is the need to cover the situation  where $\alpha(\vec{v})\sim\beta(\vec{v})$. In this case  also  $\alpha(\vec{t})\sim\beta(\vec{t})$, and the choice from  $\{\alpha(\vec{v}), \beta(\vec{v})\}$, as well as from $\{\alpha(\vec{t}), \beta(\vec{t})\}$, is indifferent. So if $f(\alpha(\vec{v}),\beta(\vec{v}))=\alpha(\vec{v})$, while  $f(\alpha(\vec{t}),\beta(\vec{t}))=\beta(\vec{t})$, then  $f(\alpha(\vec{t}),\beta(\vec{t}))\neq \alpha(\vec{t})$ while
$f(\alpha(\vec{t}),\beta(\vec{t}))\sim \alpha(\vec{t})$.

\vskip 0.2in

Unfortunately,  no choice function $f:[Fml(L)]^2\rightarrow Fml(L)$ can be uniform, for any first-order language $L$, so the definition \ref{D:uniff} is void.

\begin{Prop} \label{F:incon}
For any language  $L$ there is no uniform choice function  for $L$.
\end{Prop}

{\em Proof.} Let $L$ be any first-order language. Clearly we can pick a formula $\alpha(v)$ and variables $v_1,v_2$ such that $\alpha(v_1)\not\sim \alpha(v_2)$. Suppose  $f$ is a uniform choice function for $[Fml(L)]^2$, that is  $f$ satisfies (\ref{E:crucial}). In particular this holds for the pair $\{\alpha(v_1),\alpha(v_2)\}$ and the tuple of terms $\vec{t}=\langle v_2,v_1\rangle$.   Assume without loss of generality that  $f(\alpha(v_1),\alpha(v_2))=\alpha(v_1)$. Then  $$[f(\alpha(v_1),\alpha(v_2))](\vec{t})=\alpha(v_1)(\vec{t})=\alpha(v_2).$$
By (\ref{E:crucial})
$$[f(\alpha(v_1),\alpha(v_2))](\vec{t})\sim f(\alpha(v_1)(\vec{t}), \alpha(v_2)(\vec{t})) =f(\alpha(v_2),\alpha(v_1)),$$
therefore, by the above relations
\begin{equation} \label{E:t1}
f(\alpha(v_2),\alpha(v_1))\sim \alpha(v_2).
\end{equation}
But $f(\alpha(v_2),\alpha(v_1))=f(\alpha(v_1),\alpha(v_2))=\alpha(v_1)$ by our assumption. So $\alpha(v_1)\sim \alpha(v_2)$,  a contradiction.  \telos

\section{The sentence choice semantics (SCS) for first-order superposition logic}

We come now to examine the semantics SCS for FOLS  based on the truth relation $\models_s^2$ roughly described as Option 2 in section 2. As already said there, this semantics presumes that a restriction is imposed to the  syntax of $L_s$, namely that $|$ should not apply to quantified formulas, unless they are classical. So below we shall deal with a class of formulas of $L_s$, called ``restricted formulas/sentences''. To define them we define first the class of ``basic formulas/sentences''.

\begin{Def} \label{D:restricted}
{\em Let $L$ be a first-order language and  let ${\cal M}=\langle M,\ldots\rangle$ be an $L$-structure.

(i) The set $BFml(L_s(M))$ of} basic  formulas {\em of $L_s(M)$ is the smallest set of formulas $X$ such that (a) $Fml(L(M))\subset X$ and (b)  $X$ is closed with respect to the connectives $\wedge$, $\vee$, $\rightarrow$, $\leftrightarrow$, $|$ and $\neg$ (but not with respect to quantifiers). The set $BSen(L_s(M))$ of} basic  sentences {\em of $L_s(M)$ is the subset of $BFml(L_s(M))$ of formulas without free variables.

(ii) The set $RFml(L_s(M))$ of} restricted formulas {\em  of $L_s(M)$ is the smallest set of formulas $X$ such that (a) $BFml(L_s(M))\subset X$ and (b)  $X$ is closed with respect to $\wedge$, $\vee$, $\rightarrow$, $\leftrightarrow$,  $\neg$ and $\forall$, $\exists$ (but not with respect to $|$). The set $RSen(L_s(M))$ of} restricted  sentences {\em of $L_s(M)$ is the subset of $RFml(L_s(M))$ of formulas without free variables.}
\end{Def}

We come next to choice functions. In contrast to the choice functions used in  FCS, the choice functions $f$ of SCS apply only to pairs of  sentences of $L(M)$, i.e., $f:[Sen(L(M))]^2\rightarrow Sen(L(M))$. Let us denote  by ${\cal F}_M$ the class of all these functions. Further, SCS differs from FCS in that the collapsing function $\overline{f}$ induced from $f$ will be defined only for {\em basic sentences} i.e., for elements of $BSen(L_s(M))$.

\begin{Def} \label{D:basic}
{\em Given $f:[Sen(L(M))]^2\rightarrow Sen(L(M))$, the function $\overline{f}:BSen(L_s(M))\rightarrow Sen(L(M))$ is defined along the clauses of  Definition \ref{D:collapse} as follows:

(i) $\overline{f}(\alpha)=\alpha$, for $\alpha\in Sen(L(M))$,

(ii) $\overline{f}(\varphi\wedge\psi)=\overline{f}(\varphi)\wedge\overline{f}(\psi)$,

(iii) $\overline{f}(\neg\varphi)=\neg\overline{f}(\varphi)$,

(iv) $\overline{f}(\varphi|\psi)=f(\overline{f}(\varphi),\overline{f}(\psi))$.}
\end{Def}

\noindent It is easy to check that this definition is good, and yields $\overline{f}(\varphi)$ for every $\varphi\in BSen(L_s(M)$. Especially concerning step  (iv), note that  if $\varphi|\psi$  belongs to $BSen(L_s(M))$, then so do $\varphi$ and $\psi$,  hence $\overline{f}(\varphi)$ and $\overline{f}(\psi)$ are defined and are sentences of $L(M)$. Therefore $f(\overline{f}(\varphi),\overline{f}(\psi))$ is defined too.

We come to the truth definition $\langle {\cal M},f\rangle \models_s^2\varphi$, where ${\cal M}$ is an $L$-structure, $f\in {\cal F}_M$ and $\varphi\in RSen(L_s(M))$.

\begin{Def} \label{D:rtruth}
{\em $\langle {\cal M},f\rangle \models_s^2\varphi$ is defined
by induction on the length of $\varphi$ along the following clauses. (We think of $\wedge$, $\neg$, $|$ and $\forall$ as basic connectives, the others being thought of as abbreviations.)

(i) $\langle {\cal M},f\rangle \models_s^2\alpha$ iff ${\cal M}\models\alpha$, for $\alpha\in Sen(L(M))$.

(ii) $\langle {\cal M},f\rangle \models_s^2\varphi\wedge\psi$ iff $\langle {\cal M},f\rangle \models_s^2\varphi$ and $\langle {\cal M},f\rangle \models_s^2\psi$.

(iii) $\langle {\cal M},f\rangle \models_s^2\neg \varphi$ iff $\langle {\cal M},f\rangle \not\models_s^2\varphi$.

(iv) $\langle {\cal M},f\rangle \models_s^2\varphi|\psi$ iff ${\cal M} \models f(\overline{f}(\varphi),\overline{f}(\psi))$.

(v) $\langle {\cal M},f\rangle \models_s^2(\forall v)\varphi(v)$ iff $\langle {\cal M},f\rangle \models_s^2\varphi(x)$ for every $x\in M$.}
\end{Def}

It is easy to check that the above definition assigns a unique truth value to every $\varphi\in RSen(L_s(M))$. Specifically  clauses (i)-(iv)  attribute truth values to all basic sentences, while clause (v) is needed for quantified (non-classical) sentences.

\vskip 0.1in

Given a class $X\subseteq {\cal F}$, the  $X$-logical consequence relation $\models_X^2$ and the notion of  $X$-tautology, $\models_X^2\varphi$, are defined as usual. We denote by $Taut_X^2$ the set of $\models_X^2\varphi$-tautologies. We denote again by $Asso$, $Reg$, $Reg^*$ and $Dec$ the classes of associative, regular, regular and associative, and regular, associative and $\neg$-decreasing elements of ${\cal F}$. In particular we have
$$Dec\subset Reg^*\subset Reg\subset {\cal F},$$
hence
$$Taut({\cal F})\subseteq  Taut(Reg)\subseteq Taut(Reg^*)\subseteq Taut(Dec).$$
Given a class $X\subseteq {\cal F}$ and a formal system  $\Lambda$ consisting of axioms (set $\models_X^2\varphi$-tautologies) and rules of inference, we shall denote by
$${\rm RFOLS}(X,\Lambda)$$
the logical system having as usual semantic part $X$ and syntactic part $\Lambda$ (the prefix ``R'' stands for reminding that we work in a restricted class of sentences of $L_s$).

\subsection{Soundness}
The formal systems $\Lambda_0$, $\Lambda_1$,  $\Lambda_2$ and  $\Lambda_3$ described in section 2.1,  are going  to formalize the classes ${\cal F}$, $Reg$, $Reg^*$  and  $Dec$ of choice functions. So we shall be dealing with the logics
$${\rm RFOLS}({\cal F},\Lambda_0), \  {\rm RFOLS}(Reg,\Lambda_1), \  {\rm RFOLS}(Reg^*,\Lambda_2), \ {\rm RFOLS}(Dec,\Lambda_3).$$
Since we work with restricted formulas only we must be careful with the syntax of the systems $\Lambda_i$ above. Namely the following remarks are in order.

1) The formulas that can be substituted  in the axiom schemes $S_i$  above must be  restricted, that is $S_i\subset RFml(L_s)$.

2) Whenever  we  write $\Sigma\vdash_{\Lambda_i}\varphi$, it is implicitly assumed that $\Sigma\cup\{\varphi\}\subset RFml(L_s)$

3) Next the rule $SV$ says that if $\varphi\leftrightarrow\psi$ is provable in $\Lambda_0$, then we can derive  that $\varphi|\sigma\leftrightarrow\psi|\sigma$. But since $|$ applies only to non-quantified formulas (unless they are classical), $\varphi$ and $\psi$, hence also $\varphi\leftrightarrow\psi$, and $\sigma$ must be basic formulas.

4) Given that the above conditions are satisfied,  if $\varphi_1,\ldots,\varphi_n$ is a $\Lambda_i$-proof of $\varphi$ from $\Sigma$, then every $\varphi_i$ is restricted.

Most of the results given in this and later subsections have proofs similar  to proofs of corresponding results of \cite{Tz18}. However the adaptations needed, especially in the proofs of completeness theorems,  are rather extensive and so we give them here in full detail.

\begin{Thm} \label{T:sound}
Let $X\subseteq {\cal F}$. If  $\Lambda$ is a system such that   $\textsf{Ax}(\Lambda)\subset Taut(X)$ and $\textsf{IR}(\Lambda)=\{\textit{MP,GR}\}$, then ${\rm RFOLS}(X,\Lambda)$ is sound. In particular ${\rm RFOLS}({\cal F},\Lambda_0)$ is sound.
\end{Thm}

{\em Proof.} Let $X$, $\Lambda$ be as stated and let $\Sigma\vdash_\Lambda\varphi$, for a set of sentences $\Sigma$ and a sentence $\varphi$.
Let $\varphi_1,\ldots,\varphi_n$, where $\varphi_n=\varphi$,  be a $\Lambda$-proof of $\varphi$. As usual we show that  $\Sigma\models_X\varphi_i$, for every $1\leq i\leq n$,  by induction on $i$. Given $i$, suppose the claim holds for all $j<i$, and let $\langle {\cal M},f\rangle\models^2_s\Sigma$, for some $L$-structure ${\cal M}$ and $f\in X$. We show that  $\langle {\cal M},f\rangle\models^2_s\varphi_i$. If $\varphi_i\in \Sigma$ this is obvious. If $\varphi_i\in \textsf{Ax}(\Lambda)$, then $\langle {\cal M},f\rangle\models^2_s\varphi_i$, because by assumption $\textsf{Ax}(\Lambda)\subset Taut(X)$ and $f\in X$. Next suppose $\varphi_i$ is derived by the help of $MP$. Then there are sentences $\varphi_j$, $\varphi_k=(\varphi_j\rightarrow \varphi_i)$, for some $j,k<i$. By the induction assumption,  $\langle{\cal M},f\rangle\models^2_s\varphi_j$ and $\langle {\cal M},f\rangle\models^2_s\varphi_k$. Therefore  $\langle {\cal M},f\rangle\models^2_s\varphi_i$. Finally let $\varphi_i$ be derived by the help of $GR$, i.e., there is $j<i$ and $\varphi_j(v)$ such that $\varphi_i=(\forall v)\varphi_j(v)$. ($\Sigma$ is a set of sentences so  $v$ does not occur free in $\Sigma$.) By the induction assumption $\langle{\cal M},f\rangle\models^2_s\varphi_j(x)$ for every $x\in M$.  Then by the definition of $\models_s^2$, $\langle{\cal M},f\rangle\models^2_s(\forall v)\varphi_j(v)$. Therefore  $\langle {\cal M},f\rangle\models^2_s\varphi_i$. \telos

\vskip 0.2in

In contrast to $\Lambda_0$, the formal systems $\Lambda_i$ for $i=1,2,3$ contain in addition the rule $SV$, already mentioned in sections 1.1 and 2.1. Since however we are working in a language with syntactic restrictions  we must specify it even more concretely. Recall that $\varphi|\psi$ makes sense only if $\varphi$ and $\psi$ are basic formulas, so $SV$ takes here the form:
$$(SV) \ \ \mbox{\em For} \ \ \varphi,\psi,\sigma\in BFml(L_s), \ \mbox{if} \ \ \varphi\leftrightarrow\psi \ \mbox{is provable in $\Lambda_0$}$$
$$\mbox{\em infer that}\quad  \varphi|\sigma\leftrightarrow\psi|\sigma.$$

\begin{Thm} \label{T:soundstrong}
Let $X\subseteq Reg$. If  $\Lambda$ is a system such that   $\textsf{Ax}(\Lambda)\subset Taut(X)$ and $\textsf{IR}(\Lambda)=\{\textit{MP},GR,SV\}$, then ${\rm RFOLS}(X,\Lambda)$ is sound. In particular \newline ${\rm RFOLS}(Reg,\Lambda_1)$,
${\rm RFOLS}(Reg^*,\Lambda_2)$ and ${\rm RFOLS}(Dec,\Lambda_3)$ are sound.
\end{Thm}

{\em Proof.}  Let  $X\subseteq Reg$, $\textsf{Ax}(\Lambda)\subset Taut(X)$ and $\textsf{IR}(\Lambda)=\{\textit{MP},GR,SV\}$, and let $\Sigma\vdash_\Lambda\varphi$. Let $\varphi_1,\ldots,\varphi_n$, where $\varphi_n=\varphi$,  be a $\Lambda$-proof of $\varphi$. We show, by induction on $i$,  that for all $i=1,\ldots,n$,  $\Sigma\models_X\varphi_i$. Let $\langle {\cal M},f\rangle\models^2_s\Sigma$, with $f\in X$.
Given $\varphi_i$, the proof that $\langle {\cal M},f\rangle\models^2_s\varphi_i$ (given the induction assumption) goes exactly as in the proof of Theorem \ref{T:sound}, except of the case where $\varphi_i$ follows from a sentence $\varphi_j$, for $j<i$, by the rule $SV$. It means that $\varphi_i=(\sigma|\tau\leftrightarrow \rho|\tau)$ while $\varphi_j=(\sigma\leftrightarrow \rho)$, where $\vdash_{\Lambda_0}(\sigma\leftrightarrow \rho)$. Moreover $\sigma$, $\rho$ and $\tau$ are basic sentences. Now $\Lambda_0$ is a system satisfying the conditions of \ref{T:sound} above for $X={\cal F}$,  so $\models^2_{\cal F}(\sigma\leftrightarrow \rho)$. It means that for every $L$-structure ${\cal N}$ and every $g\in {\cal F}$, $\langle {\cal N},g\rangle\models^2_s(\sigma\leftrightarrow \rho)$. Since $\sigma$, $\rho$ and $\tau$ are basic sentences, $\overline{f}(\sigma)$, $\overline{f}(\rho)$ and $\overline{f}(\tau)$ are defined and moreover $\langle {\cal N},g\rangle\models^2_s(\sigma\leftrightarrow \rho)$ is equivalent to ${\cal N}\models\overline{g}(\sigma)\leftrightarrow \overline{g}(\rho)$. Since this holds for every ${\cal N}$,   $\overline{g}(\sigma)\leftrightarrow \overline{g}(\rho)$ is a classical tautology, or $\overline{g}(\sigma)\sim\overline{g}(\rho)$, for every $g\in {\cal F}$. In particular, $\overline{f}(\sigma)\sim\overline{f}(\rho)$. Now since $X\subseteq Reg$, $f\in X$ implies   $f$ is regular. Therefore   $\overline{f}(\sigma)\sim\overline{f}(\rho)$ implies that  $f(\overline{f}(\sigma),\overline{f}(\tau))\sim f(\overline{f}(\rho),\overline{f}(\tau))$, or $\overline{f}(\sigma|\tau)\sim \overline{f}(\rho|\tau)$, therefore ${\cal M}\models \overline{f}(\sigma|\tau)\leftrightarrow \overline{f}(\rho|\tau)$, or $\langle {\cal M},f\rangle\models^2_s(\sigma|\tau\leftrightarrow \rho|\tau)$, i.e.,  $\langle {\cal M},f\rangle\models^2_s\varphi_i$, as required. The other claim follows from the fact that the logics in question clearly satisfy the criteria of the general statement. This completes the proof. \telos

\subsection{Completeness of the logic ${\rm RFOLS}({\cal F},\Lambda_0)$}

Since the rules of $\Lambda_0$ are only $MP$ and $GR$, the Deduction Theorem (DT) holds in $\Lambda_0$ so by Fact \ref{F:eqsat}, the two forms of Completeness Theorem CT1 and CT2 are equivalent,  so we can refer simply to ``completeness'' instead of CT1- or CT2-completeness. Further by the help of DT and  standard proofs, every consistent set $\Sigma$ of formulas of $L_s$ can be extended to a complete and Henkin-complete set of formulas $\Sigma^+$ in a language $L_s^+$, where $L^+\backslash L$ consists of new constants. (Recall that a set of formulas $\Sigma$ is Henkin-complete,  if whenever $\Sigma$ contains an existential formula $\exists v\varphi(v)$, then it contains also $\varphi(c)$, for some $c\in L$,   witnessing $\exists v\varphi(v)$.) Recall also that for a  consistent and complete  $\Sigma\subseteq Sen(L_s)$ the following hold: (a)  for every $\varphi\in Sen(L_s)$, $\varphi\in \Sigma$ iff $\neg\varphi\notin \Sigma$, (b)   $\varphi\wedge \psi\in\Sigma$ iff  $\varphi\in\Sigma$ and $\psi\in \Sigma$, (c) if $\Sigma\vdash_K\varphi$, then $\varphi\in \Sigma$.

Before coming to the logics introduced  in the previous subsection,  we shall give  a general criterion of satisfiability for a consistent,  complete and Henkin-complete set $\Sigma$ of sentences of $L_s$. Given such a set $\Sigma$, if we set   $\Sigma_1=\Sigma\cap Sen(L)$ (the subset of $\Sigma$ that contains the classical sentences of $\Sigma$) then obviously  $\Sigma_1$ is a consistent,  complete and Henkin-complete  set of sentences of $L$. By the Completeness Theorem of FOL, there exists an $L$-structure  ${\cal M}$ such that, for every $\alpha\in Sen(L)$, $\alpha\in \Sigma_1$ iff ${\cal M}\models\alpha$. We  have the following criterion of  satisfiability.

\begin{Lem} \label{L:Henkin}
Let $X\subseteq {\cal F}$ and $\Lambda\subset Taut(X)$. Let also $\Sigma$ be a $\Lambda$-consistent, complete and Henkin-complete set of sentences of  $L_s$ and let  $\Sigma_1=\Sigma\cap Sen(L)$   and  ${\cal M}$ such that
\begin{equation} \label{E:constant}
\alpha\in \Sigma_1 \ \Leftrightarrow \ {\cal M}\models\alpha.
\end{equation}
Then given  $f\in X$, $\langle  {\cal M},f\rangle\models_s^2\Sigma$ if and only if  for every $\varphi\in BSen(L_s)$ (the set of basic sentences of $L_s$),
\begin{equation} \label{E:kalisynthiki}
\varphi\in \Sigma \ \Rightarrow \overline{f}(\varphi)\in \Sigma.
\end{equation}
(Actually (\ref{E:kalisynthiki})  is equivalent to $$\varphi\in \Sigma \ \Leftrightarrow \overline{f}(\varphi)\in \Sigma,$$  but the other direction follows from (\ref{E:kalisynthiki}), the consistency and completeness of $\Sigma$ and the fact that $\overline{f}(\neg\varphi)=\neg\overline{f}(\varphi)$.)
\end{Lem}

{\em Proof.}  Pick an $f\in X$ and suppose  $\langle {\cal M},f\rangle\models^2_s\Sigma$. Then  by the completeness of $\Sigma$ and the definition of $\models^2_s$, for every $\varphi\in BSen(L_s)$,
$$\varphi\in \Sigma \ \Leftrightarrow \ \langle {\cal M},f\rangle\models^2_s\varphi \Leftrightarrow {\cal M}\models \overline{f}(\varphi).$$
Now by (\ref{E:constant}), ${\cal M}\models \overline{f}(\varphi)\Rightarrow \overline{f}(\varphi)\in \Sigma_1\subset\Sigma$. Therefore $\varphi\in BSen(L_s)\cap\Sigma \ \Rightarrow \overline{f}(\varphi)\in \Sigma$.  Thus  (\ref{E:kalisynthiki}) holds.

Conversely, suppose (\ref{E:kalisynthiki}) is true.  We have to show  that  $\langle {\cal M},f\rangle\models^2_s\Sigma$. Pick some $\varphi\in \Sigma$. Assume first that $\varphi$ is basic, i.e., $\varphi\in BSen(L_s)$. Then $\overline{f}(\varphi)$ is defined. By (\ref{E:kalisynthiki}) $\overline{f}(\varphi)\in \Sigma$, therefore   $\overline{f}(\varphi)\in \Sigma_1$ since $\overline{f}(\varphi)$ is classical. So by  (\ref{E:constant}) ${\cal M}\models \overline{f}(\varphi)$. This means that $\langle {\cal M},f\rangle\models^2_s\varphi$, as required.

So it remains  to show that $\langle {\cal M},f\rangle\models^2_s\varphi$, for $\varphi\in RSen(L_s)\backslash BSen(L_s)$. In this case $\varphi$ is a Boolean and quantifier combination of basic formulas. So we can prove $\langle {\cal M},f\rangle\models^2_s\varphi$ by induction on the length of its construction from basic formulas. The steps of the induction for the Boolean connectives are trivial due to the completeness of $\Sigma$. Thus it suffices to prove $\langle {\cal M},f\rangle\models^2_s\exists v\varphi(v)$, whenever  $\exists v\varphi(v)\in \Sigma$, assuming that this is true for $\varphi$. But if $\exists v\varphi(v)\in \Sigma$, then by Henkin-completeness of $\Sigma$, $\varphi(c)\in \Sigma$ for some $c\in L$. By the induction hypothesis, $\langle {\cal M},f\rangle\models^2_s\varphi(c)$. Therefore $\langle {\cal M},f\rangle\models^2_s\exists v\varphi(v)$. This completes the proof. \telos

\vskip 0.2in

We come to the completeness of  ${\rm RFOLS}({\cal F},\Lambda_0)$. The essential step of the proof is the  Lemma.  The pattern of proof is quite similar to that of the proof of Lemma 3.8 of \cite{Tz17}.

\begin{Lem} \label{L:key}
Let $\Sigma(\vec{v})\subset RFml(L_s)$ be a $\Lambda_0$-consistent, complete and Henkin-complete set of restricted formulas of $L_s$. Then $\Sigma(\vec{v})$ is ${\cal F}$-satisfiable.
\end{Lem}

{\em Proof.} Let $\Sigma(\vec{v})$ be  a $\Lambda_0$-consistent, complete and Henkin-complete set of formulas of $L_s$. Let us set $\Sigma_1(\vec{v})=\Sigma(\vec{v})\cap Fml(L)$. Clearly $\Sigma_1(\vec{v})$ is a consistent, complete and Henkin-complete set of classical $L$-formulas. By the completeness theorem of FOL, there is an $L$-structure ${\cal M}=\langle M,\ldots\rangle$ and  a tuple $\vec{a}\in M$ such that
\begin{equation} \label{E:apply}
\alpha\in \Sigma_1(\vec{a})\Leftrightarrow {\cal M}\models\alpha.
\end{equation}
Let $L(\vec{a})$ be $L$ augmented with the parameters $\vec{a}\in M$. Without loss of generality and because of Henkin-completeness,  we may assume that $L(\vec{a})=L$. Because otherwise we may take $\Sigma_1^*(\vec{a})=\{\alpha(\vec{a}):{\cal M}\models\alpha(\vec{a})\}$ instead of $\Sigma_1(\vec{a})$. Since for every element $a_i$ of the sequence $\vec{a}$, $\exists v(v=a_i)$ belongs to $\Sigma_1^*(\vec{a})$, and $a_i$ is (essentially) the unique parameter witnessing $\exists v(v=a_i)$, it follows that $a_i\in L$. Further, if we start with $\Sigma_1^*(\vec{a})$, we can easily extend it to a consistent complete and Henkin-complete set of sentences  $\Sigma^*(\vec{a})$ of $L_s(\vec{a})$ such that  $\Sigma(\vec{a})\subseteq \Sigma^*(\vec{a})$.

Now  $\Sigma_1(\vec{a})\subset \Sigma(\vec{a})$. So, applying the criterion of Lemma \ref{L:Henkin}, in order to show that $\Sigma(\vec{a})$ is satisfiable, it suffices to construct a choice function $g:[Sen(L)]^2\rightarrow Sen(L)$ such that for every $\varphi\in BSen(L)$,
\begin{equation} \label{E:apply1}
\varphi\in \Sigma(\vec{a})\Rightarrow \overline{g}(\varphi)\in\Sigma(\vec{a}).
\end{equation}
To do that we  examine for any $\varphi,\psi\in BSen(L_s)$, the possible subsets of $\Sigma(\vec{a})$ whose elements are $\varphi|\psi$, $\varphi$, $\psi$ or  their negations. These   are the following:

\vskip 0.1in

(a1) $\{\varphi|\psi, \varphi, \psi\}\subset \Sigma(\vec{a})$

(a2)  $\{\varphi|\psi, \varphi,\neg\psi\}\subset \Sigma(\vec{a})$

(a3) $\{\varphi|\psi, \neg\varphi,\psi\}\subset \Sigma(\vec{a})$

(a4) $\{\neg(\varphi|\psi),\neg\varphi,\neg\psi\}\subset \Sigma(\vec{a})$

(a5)  $\{\neg(\varphi|\psi), \varphi,\neg\psi\}\subset \Sigma(\vec{a})$

(a6) $\{\neg(\varphi|\psi), \neg\varphi,\psi\}\subset \Sigma(\vec{a})$

\vskip 0.1in

\noindent The remaining cases,

(a7) $\{\varphi|\psi, \neg\varphi, \neg\psi\}\subset\Sigma(\vec{a})$

(a8) $\{\neg(\varphi|\psi), \varphi, \psi\}\subset\Sigma(\vec{a})$

\noindent are impossible because they contradict   $\Lambda_0$-consistency and completeness of $\Sigma(\vec{a})$. Indeed,  in case (a7) we have  $\neg\varphi\wedge\neg\psi\in \Sigma(\vec{a})$. Also  $\varphi|\psi\in\Sigma(\vec{a})$, so by  $S_2$ and completeness, $\varphi\vee\psi\in \Sigma(\vec{a})$, a contradiction. In  case (a8)  $\varphi\wedge\psi\in \Sigma(\vec{a})$. Also  $\neg(\varphi|\psi)\in\Sigma(\vec{a})$, so by $S_1$ and completeness  $\neg(\varphi\wedge\psi)\in\Sigma(\vec{a})$, a contradiction.

Given a pair $\{\alpha,\beta\}$ of sentences of $L$, we say that ``$\{\alpha,\beta\}$ satisfies   (ai)'' if for $\varphi=\alpha$ and $\psi=\beta$, the corresponding case (ai) above, for $1\leq i\leq 6$, holds. We define a choice function $g$ for $L$ as follows:
\begin{equation} \label{E:basicmap}
g(\alpha,\beta)=
\left\{\begin{array}{l}
               (i) \ \alpha, \ \mbox{if $\{\alpha,\beta\}$ satisfies  (a2) or (a6)}  \\
               (ii) \ \beta, \ \mbox{if  $\{\alpha,\beta\}$ satisfies  (a3) or (a5) } \\
               (iii) \ \mbox{any of the $\alpha$, $\beta$, if $\{\alpha,\beta\}$ satisfies (a1)
                or (a4).}
            \end{array} \right.
\end{equation}

\vskip 0.1in

{\em Claim.} $\overline{g}$ satisfies the implication (\ref{E:apply1}).

\vskip 0.1in

{\em Proof of the Claim.} We prove (\ref{E:kalisynthiki})  by induction on the length of $\varphi$. For $\varphi=\alpha\in Sen(L)$, $\overline{g}(\alpha)=\alpha$, so (\ref{E:kalisynthiki}) holds trivially. Similarly the induction steps for $\wedge$ and $\neg$ follow immediately from the fact that $\overline{g}$ commutes with these connectives and the completeness of $\Sigma(\vec{a})$. So the only nontrivial step of the induction is that for $\varphi|\psi$. It suffices to assume
\begin{equation} \label{E:as1}
\varphi\in \Sigma(\vec{a})\ \Rightarrow \overline{g}(\varphi)\in\Sigma(\vec{a}),
\end{equation}
\begin{equation} \label{E:as2}
\psi\in \Sigma(\vec{a})\ \Rightarrow \overline{g}(\psi)\in\Sigma(\vec{a}),
\end{equation}
and prove
\begin{equation} \label{E:symp}
\varphi|\psi\in \Sigma(\vec{a})\ \Rightarrow \overline{g}(\varphi|\psi)\in\Sigma(\vec{a}).
\end{equation}
Assume $\varphi|\psi\in \Sigma(\vec{a})$. Then the only possible combinations of  $\varphi$, $\psi$ and their negations that can belong to $\Sigma(\vec{a})$ are those of cases (a1), (a2) and (a3) above. To prove (\ref{E:symp}) it suffices to check that $\overline{g}(\varphi|\psi)\in\Sigma(\vec{a})$ in each of these cases. Note that $\overline{g}(\varphi|\psi)=g(\overline{g}(\varphi),\overline{g}(\psi))=g(\alpha,\beta)$, where $\overline{g}(\varphi)=\alpha$ and $\overline{g}(\psi)=\beta$ are sentences of $L$, so (\ref{E:basicmap}) applies.

Case (a1):  Then  $\varphi\in \Sigma(\vec{a})$ and  $\psi\in \Sigma(\vec{a})$. By (\ref{E:as1}) and (\ref{E:as2}), $\overline{g}(\varphi)\in\Sigma(\vec{a})$ and $\overline{g}(\varphi)\in\Sigma(\vec{a})$. By definition (\ref{E:basicmap}), $\overline{g}(\varphi|\psi)=g(\overline{g}(\varphi),\overline{g}(\psi))$ can be either $\overline{g}(\varphi)$ or $\overline{g}(\psi)$. So in either case  $\overline{g}(\varphi|\psi)\in\Sigma(\vec{a})$.

Case (a2):  Then $\varphi\in \Sigma(\vec{a})$ and   $\neg\psi\in \Sigma(\vec{a})$. By (\ref{E:as1}) and (\ref{E:as2}), $\overline{g}(\varphi)\in\Sigma(\vec{a})$, $\overline{g}(\psi)\notin\Sigma(\vec{a})$. Also  by  (\ref{E:basicmap}), $\overline{g}(\varphi|\psi)=g(\overline{g}(\varphi),\overline{g}(\psi))
=\overline{g}(\varphi)$, thus $\overline{g}(\varphi|\psi)\in \Sigma(\vec{a})$.

Case (a3):  Then  $\neg\varphi\in \Sigma(\vec{a})$,  $\psi\in \Sigma(\vec{a})$. By (\ref{E:as1}) and  (\ref{E:as2}), $\overline{g}(\varphi)\notin\Sigma(\vec{a})$, $\overline{g}(\psi)\in\Sigma(\vec{a})$. By (\ref{E:basicmap}), $\overline{g}(\varphi|\psi)=g(\overline{g}(\varphi),\overline{g}(\psi))=
\overline{g}(\psi)$, thus $\overline{g}(\varphi|\psi)\in \Sigma(\vec{a})$. This completes the proof of the Claim.

\vskip 0.1in

It follows that condition (\ref{E:apply1}) is true,  so by Lemma \ref{L:Henkin}, since  ${\cal M}\models \Sigma_1(\vec{a})$ where $\Sigma_1(\vec{a})=\Sigma(\vec{a})\cap Sen(L)$,  $\langle {\cal M},g\rangle\models_s^2\Sigma(\vec{a})$, therefore $\Sigma(\vec{a})$ is ${\cal F}$-satisfiable. \telos

\begin{Thm} \label{T:MainC0}
{\rm (Completeness of ${\rm RFOLS}({\cal F},\Lambda_0$))}
Let $\Sigma(\vec{v})$ be a consistent set of restricted formulas of $L_s$. Then  $\Sigma(\vec{v})$ is  ${\cal F}$-satisfiable, i.e., there are ${\cal M}$, $f:[Sen(L)]^2\rightarrow Sen(L)$ and $\vec{a}\in M$ such that $\langle {\cal M},f\rangle\models_s^2\Sigma(\vec{a})$.
\end{Thm}

{\em Proof.}  Let $\Sigma(\vec{v})$ be a $\Lambda_0$-consistent set of formulas. Extend $\Sigma(\vec{v})$ to a $\Lambda_0$-consistent,  complete and Henkin-complete set of formulas of $L_s^+\supseteq L_s$ such that $\Sigma(\vec{v})\subseteq \Sigma^+(\vec{v})$.   By Lemma \ref{L:key},  $\Sigma^+(\vec{v})$ is ${\cal F}$-satisfiable. Therefore so is $\Sigma(\vec{v})$.  \telos

\subsection{Conditional completeness of the remaining systems}
Coming to completeness, as in the case of PLS, the presence of $SV$ makes the status of Deduction Theorem (DT) open. In turn  the absence of DT has two consequences: (a) we don't know if CT1 and CT2 are equivalent (we only know that CT1 implies CT2) and (b) we don't know if a consistent set of formulas can be extended to a consistent and complete set (and a fortiori if it can be extended to a consistent, complete and Henkin-complete set). So, concerning the completeness of the systems based on $\Lambda_i$, for $i=1,2,3$, (a)  we shall be confined  to the weaker form CT2 only, and (b)  we shall appeal to an {\em extendibility principle} for the formal systems $\Lambda_i$, already used in \cite{Tz18}.

$$(cHext(\Lambda))
\hspace{.5\columnwidth minus .5\columnwidth} \ \ \mbox{\em Every $\Lambda$-consistent set of formulas of $L_s$ can be extended}\hspace{.5\columnwidth minus .5\columnwidth} \llap{}$$
\hspace{0.8in} $\mbox{\em to a $\Lambda$-consistent, complete and Henkin-complete set}$.

\vskip 0.1in

\noindent We can see  $cHext(\Lambda)$ as the conjunction of $cext(\Lambda)$ and $Hext(\Lambda)$, where $cext(\Lambda)$ says that every $\Lambda$-consistent set can be extended  to a complete $\Lambda$-consistent set, and  $Hext(\Lambda)$ says that every $\Lambda$-consistent set can be extended  to a Henkin-complete $\Lambda$-consistent set.

The following Lemma will be essential for the completeness of the aforementioned logics, proved in the next section.

\begin{Lem} \label{L:compessential}
If $\Sigma\subset Sen(L_s)$ is closed with respect to $\vdash_{\Lambda_i}$, for some $i=1,2,3$,  and $\alpha,\alpha'$ are  formulas of $L$ such that $\alpha\sim\alpha'$, then for every $\beta$, $(\alpha|\beta\leftrightarrow \alpha'|\beta)\in\Sigma$.
\end{Lem}

{\em Proof.} Let $\alpha\sim\alpha'$. Then $\vdash_{\rm FOL}\alpha\leftrightarrow\alpha'$, hence also  $\vdash_{\Lambda_0}\alpha\leftrightarrow\alpha'$. By $SV$ it follows that for every $\beta$, $\vdash_{\Lambda_i}\alpha|\beta\leftrightarrow\alpha'|\beta$. Therefore  $(\alpha|\beta\leftrightarrow \alpha'|\beta)\in\Sigma$ since $\Sigma$ is $\vdash_{\Lambda_i}$-closed. \telos

\vskip 0.2in

The following theorem is the analogue of Theorem 3.16 of \cite{Tz17}, as well as part of Theorem 4.9 of \cite{Tz18}.

\begin{Thm} \label{T:con-com-fol}
{\rm (Conditional CT2-completeness of ${\rm RFOLS}(Reg,\Lambda_1)$}
The logic ${\rm RFOLS}(Reg,\Lambda_1$) is {\rm CT2}-complete if and only if $cHext(\Lambda_1)$ is true.
\end{Thm}

{\em Proof.}  We prove first the easy direction. Assume $cHext(\Lambda_1)$ is false. Then either $cext(\Lambda_1)$ fails or $Hext(\Lambda_1)$ fails.

Assume first that $cext(\Lambda_1)$ fails. It follows that  there is a maximal $\Lambda_1$-consistent set of formulas $\Sigma(\vec{v})$ not extendible to a $\Lambda_1$-consistent and complete set. It means that there is a formula $\varphi(\vec{v})$ such that both $\Sigma(\vec{v})\cup\{\varphi(\vec{v})\}$ and $\Sigma(\vec{v})\cup\{\neg\varphi(\vec{v})\}$ are  $\Lambda_1$-inconsistent and hence unsatisfiable. Then clearly $\Sigma(\vec{v})$ cannot be satisfiable in any structure ${\cal M}$, for then  ${\cal M}$ would also satisfy $\varphi(\vec{v})$ or $\neg\varphi(\vec{v})$. Thus CT2-completeness  fails.

Next assume that  $Hext(\Lambda_1)$ fails. It follows  that  there is a maximal $\Lambda_1$-consistent set of formulas $\Sigma(\vec{v})$ not extendible to a $\Lambda_1$-consistent and Henkin-complete set. It means that $\Sigma(\vec{v})$ contains a formula $\exists u\varphi(u,\vec{v})$ such that  for any new constant $c$, $\Sigma(\vec{v})\cup\{\varphi(c,\overline{v})\}$ is ${\Lambda_1}$-inconsistent. But then $\Sigma(\vec{v})$ cannot be  satisfiable. For if ${\cal M}$ satisfies $\Sigma(\vec{v})$, then in particular $\exists u\varphi(u,\vec{v})$ is satisfied in ${\cal M}$,  so also  $\varphi(c,\vec{v})$ is satisfied in ${\cal M}$ for some $c\in M$. Therefore $\Sigma(\vec{v})\cup\{\varphi(c,\vec{v})\}$ is satisfiable, contrary to the fact that $\Sigma(\vec{v})\cup\{\varphi(c,\vec{v})\}$ is inconsistent. Thus indeed the ${\Lambda_1}$-consistent set $\Sigma(\vec{v})$ is not  satisfiable, so CT2-completeness  fails.

We come to the main direction of the equivalence assuming $cHext(\Lambda_1)$ is true. Then given a $\Lambda_1$-consistent set $\Sigma(\vec{v})\subset RFml(L_s)$ of restricted formulas, we may assume without loss of generality that it is also complete and Henkin-complete.  We have to find ${\cal M}$ and $g\in Reg$ such that $\langle {\cal M},g\rangle\models^2_s\Sigma(\vec{a})$ for some $\vec{a}\in M$. It  turns out that the main argument of Lemma \ref{L:key}, concerning the definition of the choice  function $g$, works  also here, with the necessary adjustments.
Namely it suffices  to find  a  choice function $g\in Reg$  such that  $\langle {\cal M},g\rangle\models_s^2\Sigma(\vec{a})$, where ${\cal M}$ and $\vec{a}\in M$ are the model and parameters such that for every $\alpha\in L(\vec{v})$,
$$\alpha\in \Sigma_1(\vec{a})\Leftrightarrow {\cal M}\models\alpha,$$
where $\Sigma_1(\vec{a})=\Sigma(\vec{a})\cap Sen(L)$. The definition of  $g$ follows exactly the pattern of definition of $g$ in the proof of Lemma \ref{L:key}, except that we need now to take care  so that $g$ be regular. Recall that $g$ is regular if for all $\alpha$, $\alpha'$, $\beta$, $$\alpha'\sim\alpha \ \Rightarrow \ g(\alpha',\beta)\sim g(\alpha,\beta).$$ In (\ref{E:basicmap})  $g$ is defined by three clauses: (i) (a2) or (a6), (ii) (a3) or (a5), (iii) (a1) or (a4).

\vskip 0.1in

{\em Claim. } The regularity constraint is satisfied whenever $g$ is defined by clauses (i) and (ii) above.

\vskip 0.1in

{\em Proof of Claim.} Pick $\alpha$, $\alpha'$, $\beta$ such that $\alpha\sim\alpha'$. We prove the Claim for the case that  $g(\alpha,\beta)$ is defined according to  clause (i)-(a2). All other cases are verified  similarly.  That $g(\alpha,\beta)$  is defined by case (i)-(a2) of (\ref{E:basicmap})  means that $\alpha|\beta\in \Sigma(\vec{a})$, $\alpha\in \Sigma(\vec{a})$, $\neg\beta\in \Sigma(\vec{a})$ and $g(\alpha,\beta)=\alpha$.  It suffices to see that necessarily $g(\alpha',\beta)=\alpha'\sim g(\alpha,\beta)$.

Since $\Sigma(\vec{a})$ is complete, it is closed with respect to $\vdash_{\Lambda_1}$, so by Lemma \ref{L:compessential}, $\alpha\sim\alpha'$ implies that $(\alpha|\beta\leftrightarrow \alpha'|\beta)\in\Sigma(\vec{a})$.  Also by assumption, $\alpha|\beta\in \Sigma(\vec{a})$, hence $\alpha'|\beta\in \Sigma(\vec{a})$. Moreover $\alpha'\in \Sigma(\vec{a})$, since $\alpha\in \Sigma(\vec{a})$, and $\neg\beta\in\Sigma(\vec{a})$. Therefore case (i)-(a2) occurs too for $\alpha'|\beta$, $\alpha'$ and $\beta$. So, by (\ref{E:basicmap}),  $g(\alpha',\beta)=\alpha'$, therefore $g(\alpha',\beta)\sim g(\alpha,\beta)$. This proves the Claim.

\vskip 0.1in

It follows from the Claim that if we define  $g$ according to  (\ref{E:basicmap}), regularity  is  guaranteed unless   $g(\alpha,\beta)$ is given by  clause (iii), that is, unless  (a1) or (a4) is the case. In such a case either both $\alpha$, $\beta$ belong to $\Sigma$, or both $\neg\alpha$, $\neg\beta$ belong to $\Sigma$, and   (\ref{E:basicmap}) allows $g(\alpha,\beta)$ to be {\em any}  of the elements $\alpha$, $\beta$. So at this point  we must intervene by a new condition that will  guarantee regularity. This is done as follows.

Pick from each $\sim$-equivalence class $[\alpha]$, a representative $\xi_\alpha\in[\alpha]$. Recall that, by completeness, the set $\Sigma_1=\Sigma\cap Sen(L)$ as well as its complement $\Sigma_2=Sen(L)-\Sigma_1$ are saturated with respect to $\sim$, that is, for every $\alpha$, either $[\alpha]\subset \Sigma_1$ or $[\alpha]\subset \Sigma_2$. Let $D_1=\{\xi_\alpha:\alpha\in \Sigma_1\}$, $D_2=\{\xi_\alpha:\alpha\in \Sigma_2\}$. Let $[D_i]^2$ be the set of pairs of elements  of $D_i$, for $i=1,2$, and pick  an arbitrary choice function $g_0:[D_1]^2\cup [D_2]^2\rightarrow D_1\cup D_2$.
Then it suffices to define $g$ by slightly revising definition  (\ref{E:basicmap}) as follows:

\begin{equation} \label{E:basicmap1}
g(\alpha,\beta)=
\left\{\begin{array}{l}
               (i) \ \alpha, \ \mbox{if $\{\alpha,\beta\}$,  satisfies  (a2) or (a6)}  \\
               (ii) \ \beta, \ \mbox{if  $\{\alpha,\beta\}$ satisfies  (a3) or (a5) } \\
               (iii) \ \sim g_0(\xi_{\alpha},\xi_{\beta}), \  \mbox{if $\{\alpha,\beta\}$ satisfies (a1) or (a4).}
            \end{array} \right.
\end{equation}
(The third clause is  just a shorthand for: $g(\alpha,\beta)=\alpha$ if $g_0(\xi_{\alpha},\xi_{\beta})=\xi_{\alpha}$, and $g(\alpha,\beta)=\beta$ if $g_0(\xi_{\alpha},\xi_{\beta})=
\xi_{\beta}$.)
In view of the Claim and the specific definition of $g$ by (\ref{E:basicmap1}), it follows immediately that  if $\alpha\sim \alpha'$ then for every  $\beta$, $g(\alpha,\beta)\sim g(\alpha',\beta)$.  So $g$ is regular. Further,  exactly  as in Lemma \ref{L:key} it follows that $\langle M,g\rangle\models^2_s\Sigma(\vec{a})$. This completes the proof.  \telos

\vskip 0.2in

The next two theorems are cited without proofs. They  are analogues of Theorems 3.18 and 3.19 of \cite{Tz17}, and their proofs follow the patterns of the latter with adaptations similar to the ones we used in the proofs of  Theorems \ref{T:MainC0} and \ref{T:con-com-fol} above.

\begin{Thm} \label{T:MainC2}
{\rm (Conditional CT2-completeness  for ${\rm RFOLS}(Reg^*,\Lambda_2$))} The logic ${\rm RFOLS}(Reg^*,\Lambda_2$) is {\rm CT2}-complete if and only if $cHext(\Lambda_2)$ is true.
\end{Thm}

\begin{Thm} \label{T:MainC3}
{\rm (Conditional CT2-completeness  for ${\rm RFOLS}(Dec,\Lambda_3$))} The logic ${\rm RFOLS}(Dec,\Lambda_3$) is {\rm CT2}-complete if and only if $cHext(\Lambda_3)$ is true.
\end{Thm}

We shall close this section and the paper by  answering a question raised in \cite{Tz17} (section 5 concerning  future work), namely whether the extension of PLS to FOLS might help us to pass from superposition of sentences to {\em superposition of objects.} Such a notion may sound a little bit strange, but is closely related to  ``disjunctive objects'' (more precisely  ``disjunctive multisets''), which have  already been  used  in \cite{Tz03} to provide  semantics for the Horn fragment of the multiplicative intuitionistic linear logic (ILL) augmented with additive disjunction. The following simple  every-day example motivates sufficiently the introduction of the concept. Restaurant menus  refer to entities of the form  ``steak or fish'' (upon choice), for main dish, and ``dessert or season fruit'' (upon choice and season), for exit.\footnote{In popular presentations of linear logic the first kind of disjunction is construed as ``multiplicative'' or deterministic, while the latter is construed as ``additive'' or non-deterministic.}   One can think of the term ``steak or fish'' as representing a new kind of {\em theoretical} entity, an object generated by the superposition of steak and fish.  Of course a specific customer who dines in the restaurant does not eat ``steak or fish''. They eat either  steak or fish, which are the {\em actualizations}, i.e., the  possible collapses, of the superposed object. It is true that {\em existence} of such objects seems dubious. They look unstable and temporary, since they always  collapse to their actualizations, and also  elusive since they can be handled  not in themselves, but only through their  actualizations. However, more or less,  the same is true   of  all theoretical entities: they are supposed to stand out there  elusive in  themselves for our minds, like  platonic ideas, accessible  only through their concrete physical realizations.   Notice in particular that in the case of superposed menu items, the phrase ``upon choice'' that accompanies them explicitly  indicates that our access to their physical realizations  is obtained only by the help  of a choice function.

One way to obtain (formal representation of) superposition of objects would be through the logic of superposition (namely FOLS), if in the latter one could prove that for any  two objects (constants) $a$ and $b$, there exists a {\em unique} object $c$ satisfying  the formula (in one free variable) $(v=a)|(v=b)$, i.e., if the sentence $(\forall v,u)(\exists!w)((w=v)|(w=u))$ would be a tautology. If that would be the case, we could write $c=a\!\uparrow\!\! b$ for the  unique object $c$ satisfying the  formula $(v=a)|(v=b)$, and  say that  $c$ is the  {\em superposition of $a$ and $b$}.

If the semantics FCS would not have broken down, it is easy to see that it would satisfy the above requirement, i.e., for every ${\cal M}$, for  every choice function $f$ for pairs of formulas and any $a,b\in M$ we would have $\langle {\cal M},f\rangle\models_s(\exists! v)((v=a)|(v=b))$. This is because  $\langle {\cal M},f\rangle\models_s(\exists! v)((v=a)|(v=b))$ holds iff ${\cal M}\models(\exists!v)f(v=a,v=b)$ and the latter is obviously true no matter whether $f(v=a,v=b)=(v=a)$ or $f(v=a,v=b)=(v=b)$.

However working with the semantics SCS for FOLS described in this section  we have the following situation.

\begin{Prop} \label{P:notsatisfy}
Let $L$ be a  first-order language,  ${\cal M}=\langle M,\ldots\rangle$ be an $L$-structure and $a\neq b\in M$.  Then:

(i) There are  choice functions $f\in {\cal F}$, such that   $\langle {\cal M},f\rangle\models_s^2(\exists!v)((v=a)|(v=b))$.

(ii) However there is no  $f\in Reg$  such that  $\langle {\cal M},f\rangle\models_s^2(\exists!v)((v=a)|(v=b))$.

\end{Prop}

{\em Proof.} (i)  Given $a\neq b$, clearly the only values for $v$ that might satisfy $(v=a)|(v=b)$ are $a$ or $b$. Pick $f$ such that  $f(a=a,a=b)=(a=a)$ and $f(b=b,a=b)=(a=b)$. Then clearly $\langle {\cal M},f\rangle\models_s^2(a=a)|(a=b)$, while $\langle {\cal M},f\rangle\not\models_s^2(b=a)|(b=b)$. Thus the only element of $M$ that satisfies $(v=a)|(v=b)$ in  $\langle {\cal M},f\rangle$ is $a$. Similarly, if we consider $f'$ such that $f'(a=a,a=b)=(a=b)$ and $f'(b=b,a=b)=(b=b)$, the only element of $M$ that satisfies $(v=a)|(v=b)$ in $\langle {\cal M},f'\rangle$ is $b$.

(ii) In contrast to (i), if $f\in Reg$ then, since $(a=a)\sim (b=b)$ we should have $f(a=a,a=b)\sim f(b=b,a=b)$. Therefore either
$$f(a=a,a=b)=(a=a), \ \mbox{and} \ f(b=b,a=b)=(b=b),$$ or
$$f(a=a,a=b)=f(b=b,a=b)=(a=b).$$ In the first case both $a,b$ satisfy $(v=a)|(v=b)$ in $\langle {\cal M},f\rangle$, while  in the second case none of the $a,b$, and hence no element of $M$, satisfies $(v=a)|(v=b)$. In either case $\langle {\cal M},f\rangle\not\models_s^2(\exists!v)((v=a)|(v=b))$. \telos

\vskip 0.2in

The preceding result shows that the attempt to represent superposition of objects through FOLS and its semantics SCS fails. We can only show that for any two objects $a,b$ there is an object $c$ satisfying the property $(v=a)|(v=b)$, but it is not unique. Nevertheless, superposition of objects can be introduced by an alternative way, namely through {\em mathematical}  rather than logical means. Specifically, given a first order theory $T$ in a language $L$, $T$ can be extended to a theory $T^|$ in the language $L\cup\{|\}$, where $|$ is a new binary operation (on the objects of $T$). $T^|$ (with underlying logic the usual FOL) consists of the axioms of $T$ plus some plausible axioms for $|$, analogous to the axioms $S_i$ of section 1.1, expressing idempotence, symmetry and associativity of $|$, and possibly some further properties for the objects $a|b$.

Let us note by the way that the notation $x|y$ was first used in \cite[\S 3]{Tz03}. The operation $x|y$ was defined there (for multisets and finite sets of multisets)  so that idempotence, symmetry and associativity hold, and also so that the object  $x|y$ be {\em distinct}  from both $x$ and $y$, i.e., $x|y\notin \{x,y\}$, unless $x=y$. In contrast, if $x|y$ is going to represent an entity that always collapses to either $x$ or $y$, then necessarily $x|y\in \{x,y\}$, i.e., $|$ must behave as  a choice function. Thus one can have at least two different implementations of the operation $x|y$ on objects: a ``projective'' one, such that $x|y\in \{x,y\}$, and a ``non-projective'' one, such that $x|y\notin \{x,y\}$.

\vskip 0.2in

{\bf Acknowledgements} Many thanks to two anonymous referees for several  corrections and suggestions that improved considerably the presentation of this paper.

\end{document}